\documentclass[draftclsnofoot]{IEEEtran}
\twocolumn   
\usepackage{graphicx}
\usepackage{amsmath}
\usepackage{hyperref}
\usepackage{amssymb}
\usepackage{amsthm}
\usepackage{amsopn}
\usepackage{mathrsfs}
\usepackage{mathtools}
\usepackage{stmaryrd}
\usepackage[T1]{fontenc}
\usepackage{lipsum}
\usepackage{graphicx}
\usepackage{enumerate}
\graphicspath{ {/home/user/images/} } 
\usepackage{subcaption}
\usepackage{cite}
  
\usepackage{tikz}
  \usetikzlibrary{matrix}
  \usetikzlibrary{fit}
  \usetikzlibrary{backgrounds}
  \usetikzlibrary{arrows}
  \usetikzlibrary{shapes}
  \usetikzlibrary{decorations.pathmorphing}
\usepackage{tikz-cd}

\newtheorem{lem}{Lemma}
\newtheorem{prop}{Proposition}

\begin{document}

\ifCLASSINFOpdf
\else
\fi

\title{Riemannian statistics meets random matrix theory\,: towards learning from high-dimensional covariance matrices}

\author{Salem Said, 
        Simon~Heuveline,
        and Cyrus Mostajeran
\thanks{S. Said is a CNRS research Scientist, at Laboratoire Jean Kuntzmann, Université Grenoble-Alpes}
\thanks{S. Heuveline is with the Centre for Mathematical Sciences of the University of Cambridge}
\thanks{C. Mostajeran is with the Department of Engineering of the University of Cambridge in the United Kingdom and the School of Physical and Mathematical Sciences of Nanyang Technological University (NTU) in Singapore}}

\maketitle

\begin{abstract}
 Riemannian Gaussian distributions were initially introduced as basic building blocks for learning models which aim to capture the intrinsic structure of statistical populations of positive-definite matrices (here called covariance matrices). While the potential applications of such models have attracted significant attention, a major obstacle still stands in the way of these applications\,: there seems to exist no practical method of computing the normalising factors associated with Riemannian Gaussian distributions on spaces of high-dimensional covariance matrices. The present paper shows that this missing method comes from an unexpected new connection with random matrix theory. Its main contribution is to prove that Riemannian Gaussian distributions of real, complex, or quaternion covariance matrices are equivalent to orthogonal, unitary, or symplectic log-normal matrix ensembles. This equivalence yields a highly efficient approximation of the normalising factors, in terms of a rather simple analytic expression. The error due to this approximation decreases like the inverse square of dimension. Numerical experiments are conducted which demonstrate how this new approximation can unlock the difficulties which have impeded applications to real-world datasets of high-dimensional covariance matrices. The paper then turns to Riemannian Gaussian distributions of block-Toeplitz covariance matrices. These are equivalent to yet another kind of random matrix ensembles, here called ``acosh-normal" ensembles. Orthogonal and unitary ``acosh-normal" ensembles correspond to the cases of block-Toeplitz with Toeplitz blocks, and block-Toeplitz (with general blocks) covariance matrices, respectively.  
\end{abstract}

\IEEEpeerreviewmaketitle

\section{Introduction}
The present paper aims to bring together the tools of Riemannian statistics and those of random matrix theory, in order to unlock the computational issues which have stood in the way of applying Riemannian learning models to datasets of high-dimensional covariance matrices.

Over the past few years, Riemannian Gaussian distributions, and mixtures of these distributions, were introduced as a means of modeling the intrinsic structure of statistical populations of covariance matrices~\cite{Cheng2013,Sa16,Sa17}. They have found successful applications in fields such as brain-computer interface analysis and artificial intelligence~\cite{zanini2018}\cite{mathieu2019}. However, such applications could not be pursued for covariance matrices of relatively larger dimension (for example, $50 \times 50$). Indeed, there seemed to exist no practical method of computing the normalising factors associated with Riemannian Gaussian distributions on spaces of high-dimensional covariance matrices. 

The theoretical contribution of this paper is to show that this missing method arises quite naturally, as soon as one realises that a strong connection exists between Riemannian Gaussian distributions and random matrix theory. 
Roughly speaking, Riemannian Gaussian distributions on spaces of real, complex, or quaternion covariance matrices correspond to log-normal orthogonal, unitary or symplectic matrix ensembles. These are similar to the classical, widely-known Gaussian orthogonal, unitary and symplectic ensembles, but with the normal weight function replaced with a log-normal weight function. Thanks to this new connection, the powerful tools of random matrix theory can be employed to uncover several original properties of Riemannian Gaussian distributions, especially for higher-dimensional covariance matrices. 

The present paper also extends the equivalence between Riemannian Gaussian distributions and random matrix ensembles to the case of block-Toeplitz matrices. Instead of a log-normal weight function $\exp(-\log^2(x))$, the corresponding weigh function is ``$\mathrm{acosh}$-normal" $\exp(-\mathrm{acosh}^2(x))$, where $\mathrm{acosh}$ is the inverse hyperbolic cosine.
No attempt is made to develop these additional ``$\mathrm{acosh}$-normal" ensembles, at least not for now. In fact, there are already plenty of results to derive and apply with the log-normal ensembles, which cover the important cases of real and complex covariance matrices. This is in part because log-normal ensembles (also called Sitltjes-Wigert ensembles) are amenable to a direct analytic treatment, as will be seen below.

It is also because of a rather surprising connection also noted in~\cite{SH21}. The log-normal unitary ensemble had appeared in the theoretical physics literature about twenty years ago, as a random matrix model for the Chern-Simons quantum field theory~\cite{Ma05}. A seen in the present paper, this ensemble corresponds to Riemannian Gaussian distributions on the space of complex covariance matrices. Of course, such Riemannian Gaussian distributions have nothing to do with quantum field theory, but some of their valuable properties can be obtained by carefully re-adapting already existing results, found within this theory. 

The present paper was developed independently from a theoretical physics paper, published only very recently~\cite{Ti20}. Both papers focus on the connection between Riemannian Gaussian distributions and log-normal matrix ensembles, but significantly differ in terms of their contribution and focus, as will be discussed below. The reader is also referred to \cite{Forrester1989} for previous results on asymptotic computations of log-partition functions from the random matrix theory literature, as well as the recent paper \cite{Forrester2021}, which considers global and local scaling limits for the $\beta=2$ Stieltjes-Wigert random matrix ensemble and associated physical interpretations. 

The paper is organised as follows. Section \ref{sec:background} recalls basic background material on Riemannian Gaussian distributions. Section \ref{sec:contribution} presents the main original results, obtained by applying random matrix theory to the study of Riemannian Gaussian distributions. Section \ref{sec:application} illustrates the importance of these original results to learning from datasets of high-dimensional covariance matrices. Finally, Section \ref{sec:siegel} considers Riemannian Gaussian distributions on the space of block-Toeplitz covariance matrices, and introduces the corresponding ``$\mathrm{acosh}$-normal" ensembles. Proofs of the propositions stated in Section \ref{sec:contribution} and \ref{sec:siegel} are provided in Appendix \ref{app:proofs}.

\section{Notation and background} \label{sec:background}
Let $\mathcal{P}^{\hspace{0.02cm}\scriptscriptstyle \beta}_{\scriptscriptstyle N}$ denote the space of $N \times N$ covariance matrices which are either real ($\beta = 1$), complex ($\beta = 2$) or quaternion ($\beta = 4$). Precisely, the elements of $\mathcal{P}^{\hspace{0.02cm}\scriptscriptstyle 1}_{\scriptscriptstyle N}$ are real, symmetric positive-definite matrices, while those of $\mathcal{P}^{\hspace{0.02cm}\scriptscriptstyle 2}_{\scriptscriptstyle N}$ and $\mathcal{P}^{\hspace{0.02cm}\scriptscriptstyle 4}_{\scriptscriptstyle N}$ are respectively complex and quaternion, hermitian positive-definite matrices.

Of course, $\mathcal{P}^{\hspace{0.02cm}\scriptscriptstyle \beta}_{\scriptscriptstyle N}$ is an open convex cone, sitting inside the real vector space $\mathcal{S}^{\hspace{0.02cm}\scriptscriptstyle \beta}_{\scriptscriptstyle N}$ of real, complex or quaternion self-adjoint matrices (the adjoint of a matrix being its conjugate transpose). Therefore, $\mathcal{P}^{\hspace{0.02cm}\scriptscriptstyle \beta}_{\scriptscriptstyle N}$ is a real differentiable manifold, whose tangent space $T_{\scriptscriptstyle Y}\mathcal{P}^{\hspace{0.02cm}\scriptscriptstyle \beta}_{\scriptscriptstyle N}$ at any point $Y \in \mathcal{P}^{\hspace{0.02cm}\scriptscriptstyle \beta}_{\scriptscriptstyle N}$ is isomorphic to $\mathcal{S}^{\hspace{0.02cm}\scriptscriptstyle \beta}_{\scriptscriptstyle N}\hspace{0.03cm}$. In particular, the dimension of $\mathcal{P}^{\hspace{0.02cm}\scriptscriptstyle \beta}_{\scriptscriptstyle N}$ is $\dim^{\hspace{0.02cm}\scriptscriptstyle \beta}_{\scriptscriptstyle N} = NN_\beta$ where $N_\beta = \frac{\beta}{2}(N-1) + 1$.

The elements $Y$ of $\mathcal{P}^{\hspace{0.02cm}\scriptscriptstyle \beta}_{\scriptscriptstyle N}$ are in one-to-one correspondence with the centred (zero-mean) $N$-variate normal distributions (real, complex circular, or quaternion circular, according to the value of $\beta$). Thus, $\mathcal{P}^{\hspace{0.02cm}\scriptscriptstyle \beta}_{\scriptscriptstyle N}$ can be equipped with the Rao-Fisher information metric of the centred $N$-variate normal model~\cite{amaribook}. This is identical to the so-called affine-invariant metric introduced in~\cite{Pennec2006},
\begin{equation} \label{eq:aff-inv}
\langle u\hspace{0.02cm},v\rangle_{\scriptscriptstyle Y} = \Re\,\mathrm{tr}\left(\left(Y^{\scriptscriptstyle -1}u\right)\left(Y^{\scriptscriptstyle -1}v\right)\right)
\hspace{0.4cm}
u\hspace{0.02cm},v \in T_{\scriptscriptstyle Y}\mathcal{P}^{\hspace{0.02cm}\scriptscriptstyle \beta}_{\scriptscriptstyle N} \simeq \mathcal{S}^{\hspace{0.02cm}\scriptscriptstyle \beta}_{\scriptscriptstyle N}
\end{equation}
where $\Re$ denotes the real part and $\mathrm{tr}$ denotes the trace. The Riemannian geometry of the metric (\ref{eq:aff-inv}) is quite well-known in the geometric information science community (see the recent book~\cite{Pennec2019}). Recall here the associated geodesic distance
\begin{equation} \label{eq:aff-dist}
 d^{\hspace{0.03cm} 2}(X\hspace{0.02cm},Y) = \mathrm{tr}\!\left(\log^2\!\left(X^{\scriptscriptstyle -\frac{1}{2}}Y X^{\scriptscriptstyle -\frac{1}{2}}\right)\right)
 \hspace{1cm}
X\hspace{0.02cm},Y \in \mathcal{P}^{\hspace{0.02cm}\scriptscriptstyle \beta}_{\scriptscriptstyle N}
\end{equation}
where matrix logarithms and powers are understood as self-adjoint matrix functions, obtained by taking logarithms and powers of eigenvalues. The main advantage of the metric (\ref{eq:aff-inv}) is its invariance under the action of the group $G^{\hspace{0.02cm}\scriptscriptstyle \beta}_{\scriptscriptstyle N}$ of $N \times N$ invertible matrices with real, complex or quaternion entries (according to the value of $\beta$). For example, if $A \in G^{\hspace{0.02cm}\scriptscriptstyle \beta}_{\scriptscriptstyle N}$ and $X$ is replaced by $A\cdot X = AX\!A^\dagger$ while $Y$ is replaced by $A\cdot Y = AY\!A^\dagger$, then the distance $d(X\hspace{0.02cm},Y)$ remains unchanged (note that $^\dagger$ denotes the adjoint, or conjugate transpose). 

In terms of the distance (\ref{eq:aff-dist}), we define Riemannian Gaussian distributions on the space $\mathcal{P}^{\hspace{0.02cm}\scriptscriptstyle \beta}_{\scriptscriptstyle N}$ to be given by parameterised probability density function $p(Y|\bar{Y},\sigma)$, as in~\cite{Sa16}\cite{Sa17},
\begin{equation} \label{eq:gdensity}
p(Y|\bar{Y},\sigma) = \left(Z(\sigma)\right)^{-1}\exp\left[-\frac{d^{\hspace{0.03cm} 2}(Y\hspace{0.02cm},\bar{Y})}{2\sigma^2}\right]
\end{equation}
where $\bar{Y} \in \mathcal{P}^{\hspace{0.02cm}\scriptscriptstyle \beta}_{\scriptscriptstyle N}$ is the centre-of-mass parameter, and $\sigma > 0$ the dispersion parameter. The normalising factor $Z(\sigma)$ is given by the integral
\begin{equation} \label{eq:Z1}
  Z(\sigma) = \int_{\mathcal{P}^{\hspace{0.02cm}\scriptscriptstyle \beta}_{\scriptscriptstyle N}}\,\exp\left[-\frac{d^{\hspace{0.03cm} 2}(Y\hspace{0.02cm},\bar{Y})}{2\sigma^2}\right]dv(Y)
\end{equation}
with respect to the Riemannian volume $dv(Y) = \det(Y)^{-N_\beta}\lbrace dY\rbrace$. Here, 
$$
\begin{array}{l}
\lbrace dY\rbrace_{(\beta = 1)} = \prod_{i\leq j} dY_{ij}   \\[0.1cm]
\lbrace dY\rbrace_{(\beta = 2)}  = \prod_{i\leq j} dY^{(a)}_{ij}\prod_{i < j} dY^{(b)}_{ij}   \\[0.1cm]
\lbrace dY\rbrace_{(\beta = 4)} = \prod_{i\leq j} dY^{(a)}_{ij}\prod_{i < j} dY^{(b)}_{ij}
 dY^{(c)}_{ij}
 dY^{(d)}_{ij}
\end{array}
$$
where $Y_{ij} = Y^{(a)}_{ij} + Y^{(b)}_{ij}\hspace{0.04cm}\mathrm{i}$ if $Y_{ij}$ is complex and $Y_{ij} = Y^{(a)}_{ij} + Y^{(b)}_{ij}\hspace{0.04cm}\mathrm{i} + Y^{(c)}_{ij}\hspace{0.04cm}\mathrm{j} + Y^{(d)}_{ij}\hspace{0.04cm}\mathrm{k}$ if $Y_{ij}$ is quaternion (here, $\mathrm{i},\mathrm{j},\mathrm{k}$ denote complex or quaternion imaginary units).

Both formulae (\ref{eq:gdensity}) and (\ref{eq:Z1}) are greatly simplified by introducing ``polar coordinates". Each $Y \in \mathcal{P}^{\hspace{0.02cm}\scriptscriptstyle \beta}_{\scriptscriptstyle N}$ can be diagonalised as $Y = U e^r U^\dagger$ where $e^r$ is a diagonal matrix, with diagonal elements $e^{r_i}$ for $(r_{\scriptscriptstyle 1},\ldots,r_{\scriptscriptstyle N}) \in \mathbb{R}^N$, and $U \in K^{\hspace{0.02cm}\scriptscriptstyle \beta}_{\scriptscriptstyle N\hspace{0.02cm}}$, which means $UU^\dagger = \mathrm{I}_{\scriptscriptstyle N}$ (the $N \times N$ identity matrix). Note that $K^{\hspace{0.02cm}\scriptscriptstyle 1}_{\scriptscriptstyle N\hspace{0.02cm}}$ is the orthogonal group $O(N)$, $K^{\hspace{0.02cm}\scriptscriptstyle 2}_{\scriptscriptstyle N\hspace{0.02cm}}$ is the unitary group $U(N)$, and $K^{\hspace{0.02cm}\scriptscriptstyle 4}_{\scriptscriptstyle N\hspace{0.02cm}}$ is the symplectic (quaternion unitary) group $\mathrm{Sp}(N)$ (for deeper insight, see the recent review~\cite{edelman}). 

Now, as in~\cite{Sa17}, let $Y$ follow the Riemannian Gaussian density (\ref{eq:gdensity}), and $Y^\prime = \bar{Y}^{\scriptscriptstyle -\frac{1}{2}}Y\bar{Y}^{\scriptscriptstyle -\frac{1}{2}}\hspace{0.03cm}$.
If $Y^\prime = U e^r U^\dagger$, then $U$ is uniformly distributed on $K^{\hspace{0.02cm}\scriptscriptstyle \beta}_{\scriptscriptstyle N\hspace{0.02cm}}$ (for a rigorous definition of this ``uniform distribution", see~\cite{meckes}), while $(r_{\scriptscriptstyle 1},\ldots,r_{\scriptscriptstyle N})$ have joint probability density function
\begin{equation} \label{eq:rjointpdf}
 p(r|\sigma) \propto \prod_{i} \exp\left[-\frac{r^{\scriptscriptstyle 2}_i}{2\sigma^{\scriptscriptstyle 2}}\right]\hspace{0.03cm}
 \prod_{i<j} \sinh^{\hspace{0.02cm}\beta}\left|\frac{r_i-r_j}{2}\right|
\end{equation}
where $\propto$ indicates proportionality and $\sinh$ the hyperbolic sine. In addition, the normalising factor $Z(\sigma)$ in (\ref{eq:Z1}) reduces to a certain multiple integral $z_\beta(\sigma)$. Specifically,
\begin{equation} \label{eq:Z2}
 Z(\sigma) = \Omega_N \hspace{0.03cm} z_\beta(\sigma)
\end{equation}
where $\Omega_N$ is a numerical constant, which appears after the 
uniformly distributed matrix $U$ is integrated out of (\ref{eq:Z1}), and where
\begin{equation} \label{eq:z1}
 z_\beta(\sigma) = \frac{1}{N!}\int_{\mathbb{R}^N}
 \prod_{i<j} \sinh^{\beta}\left|\frac{r_i-r_j}{2}\right|\hspace{0.02cm}\prod_{i} \exp\!\left[-\frac{r^{\scriptscriptstyle 2}_i}{2\sigma^{\scriptscriptstyle 2}}\!\right] dr_i
\end{equation}
One of the main issues addressed in the present paper is the efficient approximation of the multiple integral $z_\beta(\sigma)$. Until now (in~\cite{Sa16}\cite{Sa17}), this was done using a Monte Carlo technique which involved a smoothing method (containing certain arbitrarily fixed parameters) and which failed to produce coherent results when the dimension $N$ increased beyond $N = 20$. 



Before proceeding to present our main results, we should briefly note the existence of several alternative proposals for the extension of Gaussian distributions to Riemannian manifolds, including constructions based on heat flows and diffusion processes \cite{Sommer2015}. Further discussion of these alternative formulations is beyond the scope of this paper and the interested reader is referred to the literature on the topic for further information. See \cite{Pennec2019} (Section 3.4.3 and Chapter 10) and the references therein for a recent and comprehensive account.

\section{Main results} \label{sec:contribution}
The main results of the present paper stem from the equivalence between Riemannian Gaussian distributions on the spaces of real, complex, or quaternion covariance matrices and log-normal orthogonal, unitary, and symplectic matrix ensembles.

Let $Y$ be a random matrix in $\mathcal{P}^{\hspace{0.02cm}\scriptscriptstyle \beta}_{\scriptscriptstyle N}$ which follows the Riemannian Gaussian density $p(Y|\bar{Y},\sigma)$ of (\ref{eq:gdensity}). It is always possible to assume that $\bar{Y} = \mathrm{I}_{\scriptscriptstyle N}\hspace{0.03cm}$, since this holds after replacing $Y$ by $Y^\prime = \bar{Y}^{\scriptscriptstyle -\frac{1}{2}}Y\bar{Y}^{\scriptscriptstyle -\frac{1}{2}}\hspace{0.03cm}$. Then, the probability distribution of $Y$ is described by the following proposition.
\begin{prop} \label{prop:log-normal}
 Let $Y$ follow the Riemannian Gaussian density (\ref{eq:gdensity}) with $\bar{Y} = \mathrm{I}_{\scriptscriptstyle N}\hspace{0.03cm}$. If $X = e^{N_\beta\sigma^2}\hspace{0.02cm}Y$, then the probability distribution of $X$ is given by
 \begin{equation} \label{eq:lognormal1}
   \mathbb{P}(X\in B) \propto \int_B\,\mathrm{etr}\left[-\frac{\log^2(X)}{2\sigma^2}\right]\lbrace dX\rbrace  
 \end{equation}
 for any measurable subset $B$ of $\mathcal{P}^{\hspace{0.02cm}\scriptscriptstyle \beta}_{\scriptscriptstyle N\hspace{0.03cm}}$. Here, $\mathrm{etr}(\cdot) = \exp(\mathrm{tr}(\cdot))$ and the notation $\lbrace dX\rbrace$ was introduced after (\ref{eq:Z1}).  
\end{prop}
In plain words, this proposition states that $X$ follows a log-normal matrix ensemble. If $X$ is diagonalised as $X = Ux\hspace{0.03cm}U^\dagger$, then $U$ is uniformly distributed on $K^{\hspace{0.02cm}\scriptscriptstyle \beta}_{\scriptscriptstyle N\hspace{0.02cm}}$, and the eigenvalues $(x_{\scriptscriptstyle 1},\ldots,x_{\scriptscriptstyle N})$, which are all positive, have joint probability density function
\begin{equation} \label{eq:lognormal2}
  p(x|\sigma) \propto |V(x)|^\beta\prod_i \rho(x_i\hspace{0.02cm},2\sigma^2) 
\end{equation}
where $V(x) = \prod_{i<j}(x_j-x_i)$ is the Vandermonde determinant, and $\rho(x\hspace{0.02cm},k) = \exp(-\log^2(x)/k)$ is the log-normal weight function.

In essence, Proposition \ref{prop:log-normal} is already contained in~\cite{Ti20}. As a consequence of this proposition, the multiple integral $z_\beta(\sigma)$ of (\ref{eq:z1}) may be expressed as follows
\begin{align} \label{eq:z2}
\nonumber  z_\beta(\sigma) = \left(2\pi\sigma^2\right)^{\!\scriptscriptstyle N/2}\times\,\exp\left[-NN^2_\beta(\sigma^2/2)\right]\times \\  \frac{1}{N!}\int_{\mathbb{R}^N_+}\,|V(x)|^\beta\,\omega(dx)
\end{align}
where $\omega(dx) = \prod\, \bar{\rho}(dx_i)$ (product over $i=1,\ldots,N$) with $\bar{\rho}(dx_i) = \left(2\pi\sigma^2\right)^{\!\scriptscriptstyle -1/2}\rho(x_i\hspace{0.02cm},2\sigma^2)\hspace{0.02cm}dx_i\hspace{0.03cm}$. In~\cite{Ti20}, integrals as the one in (\ref{eq:z2}) are expressed using the Andreev or De Bruijn identities, often employed in random matrix theory. These yield somewhat cumbersome formulae, involving determinants of size $N \times N$. On the other hand, the present paper focuses on an alternative approach, which turns out to be more suitable from a practical point of view. Instead of expressing $z_\beta(\sigma)$ exactly, by means of complicated formulae, the aim is to use a highly efficient approximation, which involves a single analytic expression. This is indicated by the following proposition.
\begin{prop} \label{prop:trilog}
In the limit where $N \rightarrow \infty$ and $\sigma \rightarrow 0$, while the product $t = N\sigma^2$ remains constant,  
 \begin{align} \label{eq:trilog}
   \nonumber \frac{1}{N^2}\hspace{0.02cm}\log z_\beta(\sigma) \longrightarrow \frac{\beta}{2}\Phi\left(\frac{\beta}{2}t\right) \\[0.2cm]
   \text{where }\hspace{0.3cm}
   \Phi(\xi) = \frac{\xi}{6} - \frac{\mathrm{Li}_{\scriptscriptstyle 3}(e^{-\xi}) - \mathrm{Li}_{\scriptscriptstyle 3}(1)}{\xi^2}
 \end{align}
 Here, $\mathrm{Li}_{\scriptscriptstyle 3}$ is the trilogarithm function, $\mathrm{Li}_{\scriptscriptstyle 3}(\eta) = \sum^\infty_{k=1} \eta^k/k^3$.
\end{prop}
In~\cite{Ti20}, the limit in (\ref{eq:trilog}) is only mentioned in passing, and not stated under the same form. Here, this limit will be given centre stage. Proposition \ref{prop:trilog} states that (\ref{eq:trilog}) is valid in the ``double-scaling regime" ($N \rightarrow \infty$ and $\sigma \rightarrow 0$), but numerical experiments have shown that $\log z_\beta(\sigma)$ can be replaced by the expression afforded by (\ref{eq:trilog}) without notable loss of accuracy, whenever $\sigma$ is small in comparison with $N$ (this is further illustrated below). 

At least informally, this can be justified by appealing to arguments
originating in theoretical physics~\cite{Ma05}. Considered as a function of $t = N\sigma^2$, $F(t) = \log z_\beta(\sigma)$ is called the Free energy (log of partition function). This free energy can be expanded in an asymptotic series (see Section 1.3 of~\cite{Ma05}),
\begin{equation} \label{eq:feynman1}
  F(t) \sim \sum^\infty_{g=0} F_g(t)\hspace{0.02cm}\left(\frac{t}{N}\right)^{\!2g-2}
\end{equation}
which is obtained by summing over Feynman diagrams, or so-called fatgraphs. Each coefficient $F_g(t)$ is itself a series $F_g(t) \sim \sum_h  F_{g,h}\hspace{0.02cm}t^h$, where $F_{g,h}$ counts fatgraphs which are said to have $h$ holes and genus $g$ (this means that one thinks of a fatgraph as a graph with $h$ loops, drawn on a surface of genus $g$, such as a sphere or torus, \textit{etc}). Now, accepting (\ref{eq:feynman1}), it follows that for each fixed value of $t$,
\begin{equation} \label{eq:feynman2}
  \frac{1}{N^2}\hspace{0.02cm}F(t) = F_{\scriptscriptstyle 0}(t) + O\left(\frac{1}{N^2}\right)
\end{equation}
so that $F_{\scriptscriptstyle 0}(t)$ approximates the left-hand side up to an error of the order of $1/N^2$. Finally, recalling that $F(t) = \log z_\beta(\sigma)$, it is clear that $F_{\scriptscriptstyle 0}(t)$ is the right-hand side of (\ref{eq:trilog}) --- this follows by uniqueness of asymptotic expansions. 

In addition to the asymptotic form of $\log z_\beta(\sigma)$, another quantity of interest is the asymptotic empirical distribution of eigenvalues, of a random matrix $Y$ which follows the Riemannian Gaussian density $p(Y|\bar{Y},\sigma)$ with $\bar{Y} = \mathrm{I}_{\scriptscriptstyle N}\hspace{0.03cm}$. Let $(y_{\scriptscriptstyle 1},\ldots,y_{\scriptscriptstyle N})$ denote the eigenvalues of $Y$, and consider their empirical distribution
\begin{equation} \label{eq:R1}
  \hat{\nu}_\beta(I) = \mathbb{E}\left[\frac{\left| y_i \in I\right|}{N} \right]
\end{equation}
for any open interval $I \subset \mathbb{R}_+$, where $\mathbb{E}$ denotes expectation, and $| y_i \in I|$ the number of $y_i$ which belong to $I$. In~\cite{Ti20}, the probability density function of $\hat{\nu}_{\scriptscriptstyle 2}$ was expressed as a weighted sum of Gaussian distributions, by a direct application of the Christoffel-Darboux formula, as in~\cite{deift}. It is possible to do so for any value of $N$, only because the $\beta = 2$ case can be studied using a well-known family of orthogonal polynomials, called Stieltjes-Wigert polynomials. No such analytic tool is available when $\beta = 1$ or $4$.

To make up for this issue, the following proposition provides an asymptotic expression of the distribution $\hat{\nu}_\beta\hspace{0.03cm}$, valid for all values $\beta = 1, 2, 4$.
\begin{prop} \label{prop:empirical}
  In the limit where $N \rightarrow \infty$ and $\sigma \rightarrow 0$, while the product $t = N\sigma^2$ remains constant, the empirical distribution $\hat{\nu}_\beta$ converges weakly to a distribution with probability density function $n(y|\beta t/2)$, where
\begin{equation}\label{eq:empirical}
  n(y|\xi) = \frac{1}{\pi \xi y}\,\mathrm{arctan}\!\left(\frac{\sqrt{4e^\xi\hspace{0.02cm} y - (y+1)^2}}{y+1}\right) 
\end{equation}
on the interval $a(\xi) \leq y \leq b(\xi)$, where $a(\xi) = c(1+\sqrt{1-c})^{-2}$ and $b(\xi) = c(1-\sqrt{1-c})^{-2}$, with $c = e^{-\xi}\hspace{0.03cm}$.
\end{prop}
One hopes that, similar to the situation discussed after Proposition \ref{prop:trilog}, the asymptotic density (\ref{eq:empirical}) approximates the finite-$N$ empirical distribution $\hat{\nu}_\beta$ to such a good accuracy that one can replace $\hat{\nu}_\beta$ by this asymptotic density, for many practical purposes. This possibility will not be further investigated in the present paper.

\section{Towards learning applications} \label{sec:application}

Riemannian Gaussian distributions were initially proposed as basic building blocks for learning models which aim to capture the intrinsic structure of statistical populations of covariance matrices. These include the mixture models and hidden Markov models, introduced in~\cite{Sa16}\cite{mtnshmm} and further developed in~\cite{QT21}. In order to make use of these models in real-world applications, it is indispensable to know how to effectively compute the logarithm of the multiple integral $z_\beta(\sigma)$ of (\ref{eq:z1}),  for a dimension $N$ which may be in the tens or hundreds. 

Knowledge of $\log z_\beta(\sigma)$ is already crucial in the simplest situation, where one tries to fit a single Riemannian Gaussian density $p(Y|\bar{Y},\sigma)$ (rather than a whole mixture) to data $Y_{\scriptscriptstyle 1},\ldots, Y_{\scriptscriptstyle M} \in \mathcal{P}^{\hspace{0.02cm}\scriptscriptstyle \beta}_{\scriptscriptstyle N}$. Indeed, this will require setting $\bar{Y} = \hat{Y}_{\scriptscriptstyle M}$ and $\sigma = \hat{\sigma}_{\scriptscriptstyle M}\hspace{0.03cm}$, where $\hat{Y}_{\scriptscriptstyle M}$ and $\hat{\sigma}_{\scriptscriptstyle M}$ are the maximum-likelihood estimates, given in~\cite{Sa16}\cite{Sa17}. Specifically,   
\begin{equation} \label{eq:frechet}
 \hat{Y}_{\scriptscriptstyle M} = \mathrm{argmin}_{Y \in \mathcal{P}^{\hspace{0.02cm}\scriptscriptstyle \beta}_{\scriptscriptstyle N}}\,\sum^M_{m=1} d^{\hspace{0.03cm} 2}(Y_{\scriptscriptstyle m}\hspace{0.03cm},Y)
\end{equation}
is the Fréchet mean of the data $Y_{\scriptscriptstyle m}\hspace{0.03cm}$, with respect to the distance (\ref{eq:aff-dist}), and $\hat{\sigma}_{\scriptscriptstyle M}$ is the solution of the nonlinear equation
\begin{equation} \label{eq:siglikelihood}
  \phi(\hat{\sigma}_{\scriptscriptstyle M}) =  \frac{1}{M}\sum^M_{m=1} d^{\hspace{0.03cm} 2}(Y_{\scriptscriptstyle m}\hspace{0.03cm},\hat{Y}_{\scriptscriptstyle M}) \,;\, \phi(\sigma) = \sigma^3\hspace{0.03cm}\frac{d}{d\sigma}\log z_\beta(\sigma)
\end{equation}
Therefore, it is already impossible to solve a toy problem, with a single Riemannian Gaussian density, without having some kind of hold on $\log z_\beta(\sigma)$.

Until now (in~\cite{Sa16}\cite{Sa17}), $z_\beta(\sigma)$ was approximated using an \textit{ad hoc} Monte Carlo technique, which failed to produce coherent results for a dimension $N$ just above $N = 20$. In the present section, the aim will be to show that significantly improved results can be obtained by solving Equation (\ref{eq:siglikelihood}) after approximating $\log z_\beta(\sigma)$ with the expression afforded by (\ref{eq:trilog}), according to Proposition \ref{prop:trilog}.

First, numerical experiments were conducted to verify the validity of this new approximation. According to (\ref{eq:feynman2}), the right-hand side of (\ref{eq:trilog}) should approximate $\log z_\beta(\sigma)/N^2$ up to an error of the order of $1/N^2$. To see that this is correct, the right-hand side of (\ref{eq:trilog}) was compared to certain exact expressions of $\log z_\beta(\sigma)$. Namely, for the $\beta = 2$ case, one has the following expression, obtained in~\cite{habilitation}, 
\begin{align} \label{eq:swp}
\nonumber \log z_{\scriptscriptstyle 2}(\sigma) = \frac{N}{2} \log( 2\pi\sigma^2) + N(N^2 - 1)(\sigma^{2}/6) \,+ \\ \sum^{N-1}_{n=1}(N-n)\log\!\left(1 - e^{-n\hspace{0.02cm}\sigma^2}\right)
\end{align}
and for the $\beta = 1$ case, when the dimension $N$ is even, the following expression, based on~\cite{Ti20},
\begin{align} \label{eq:pfaff}
\nonumber \log z_{\scriptscriptstyle 1}(\sigma) = 
\frac{N}{2} \log( 2\pi\sigma^2) \,- \\ N(N+1)^2(\sigma^2/8) +  \log\mathrm{Pf}\left[M(\sigma)\right]
\end{align}
where $\mathrm{Pf}$ denotes the Pfaffian, equal to the square root of the determinant, and the matrix $M(\sigma)$ has entries
\begin{equation} \label{eq:MB1}
  M_{ij}(\sigma) = \exp\left[(i^2+j^2)(\sigma^2/2)\right]\mathrm{erf}((j-i)(\sigma/2))  
\end{equation}
for $i,j=1,\ldots,N$ and with $\mathrm{erf}$ the error function.

\begin{figure*}[t]
    \centering
    \begin{subfigure}[b]{0.475\textwidth}
        \centering
        \includegraphics[width=\textwidth]{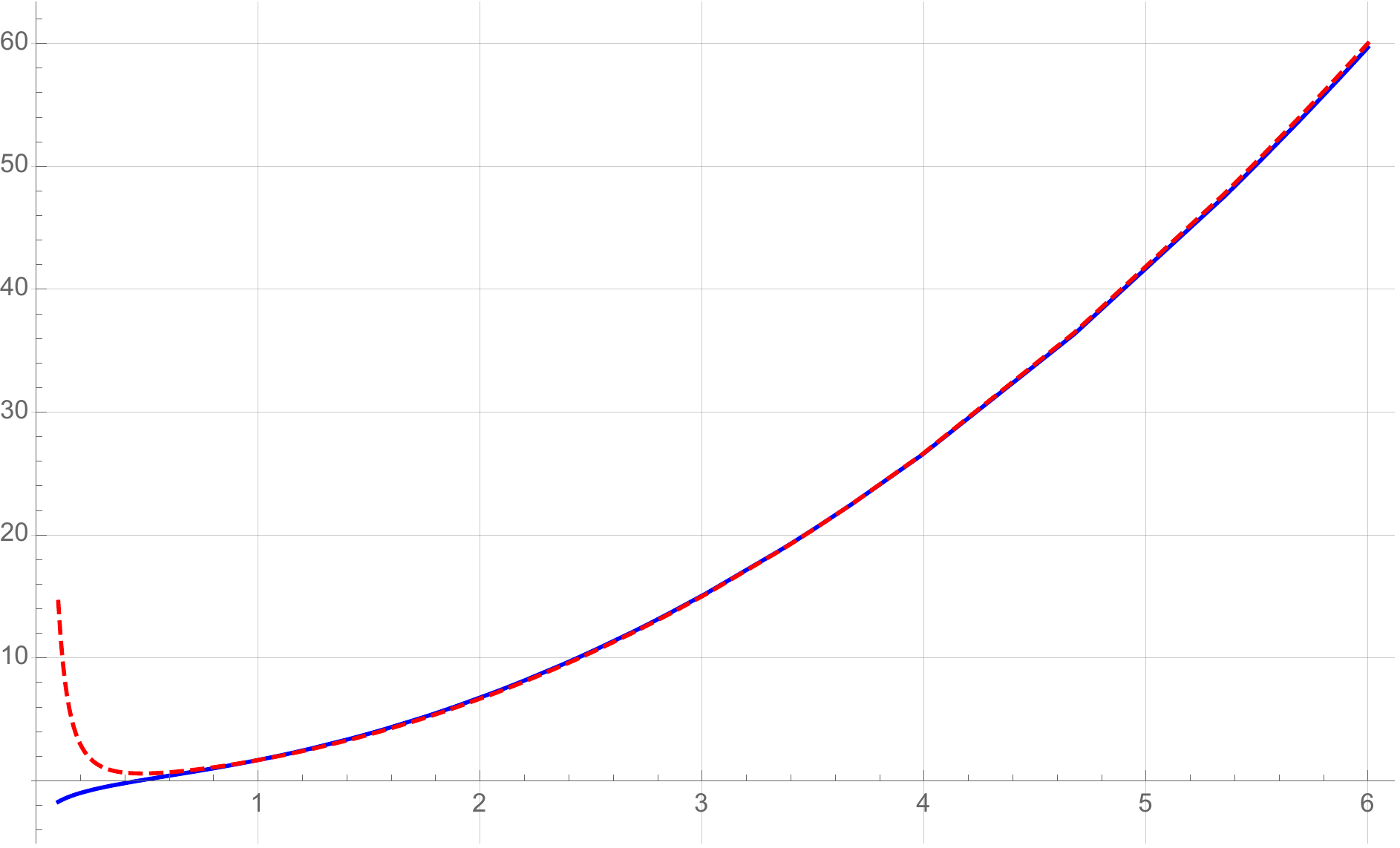}
        \caption{$\beta = 2$ and $N = 10$\,: plot of $\log z_\beta(\sigma)/N^2$}
    \end{subfigure}
    \quad
    \begin{subfigure}[b]{0.475\textwidth}  
        \centering 
        \includegraphics[width=\textwidth]{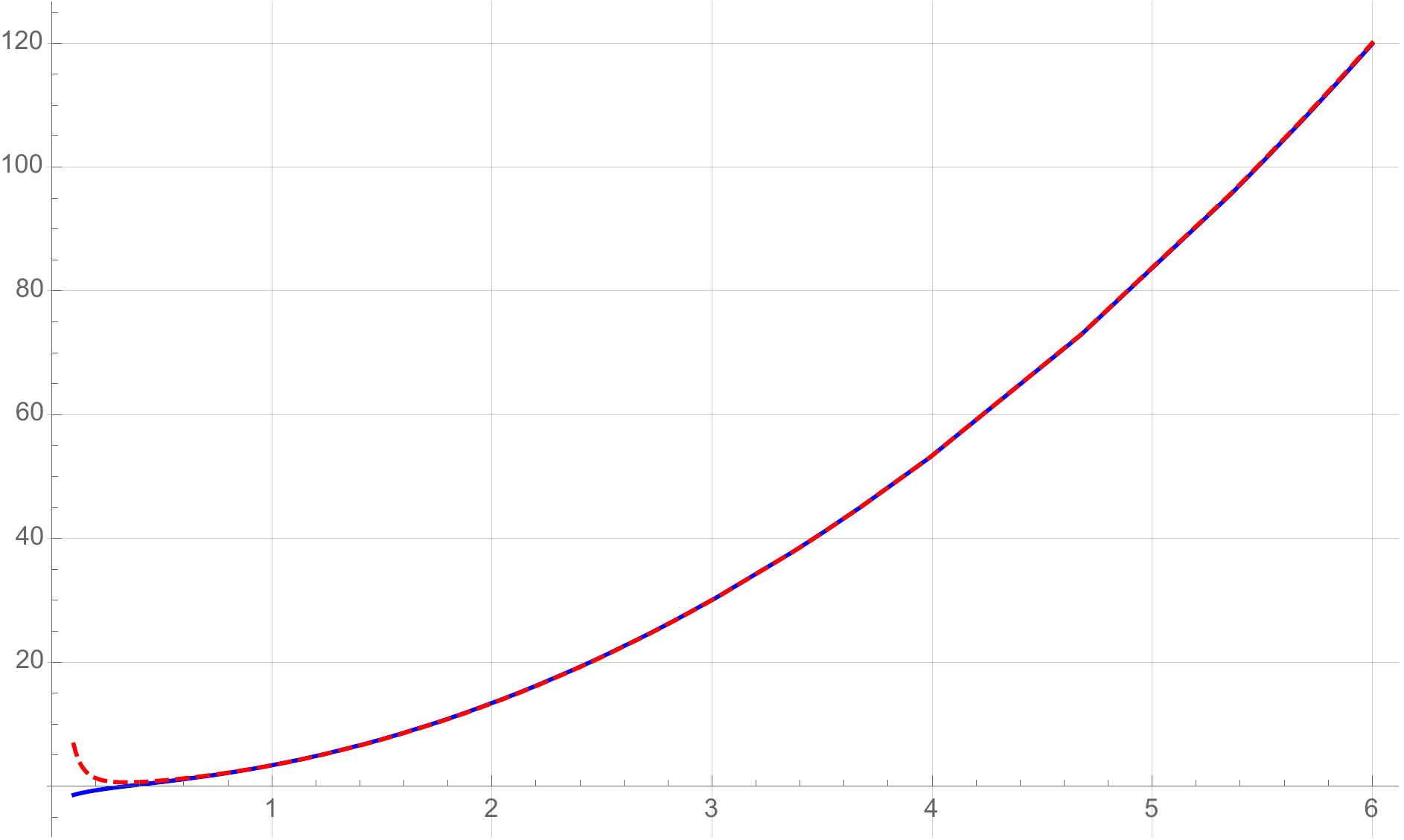}
        \caption{$\beta = 2$ and $N = 20$\,: plot of $\log z_\beta(\sigma)/N^2$}
    \end{subfigure}
    \caption{Proposition \ref{prop:trilog} ($\beta = 2$)\,: (\ref{eq:trilog}) in dashed red and (\ref{eq:swp}) in solid blue}
\label{fig:1}
\end{figure*}

Numerical evaluation of (\ref{eq:swp}) for large values of $N$ or $\sigma$ (up to $\sigma = 10$) is quite straightforward. Moreover, it immediately shows that (\ref{eq:trilog}) and (\ref{eq:swp}) agree very closely when $\sigma$ is smaller than $N$, and then gradually diverge away from one another as $\sigma$ increases. Figure \ref{fig:1} provides graphical illustration for $N = 10$ and $20$. Still larger values of $N$ yield an even stronger match between (\ref{eq:trilog}) and (\ref{eq:swp}).

\begin{figure*}[t]
    \centering
    \begin{subfigure}[b]{0.475\textwidth}
        \centering
        \includegraphics[width=\textwidth]{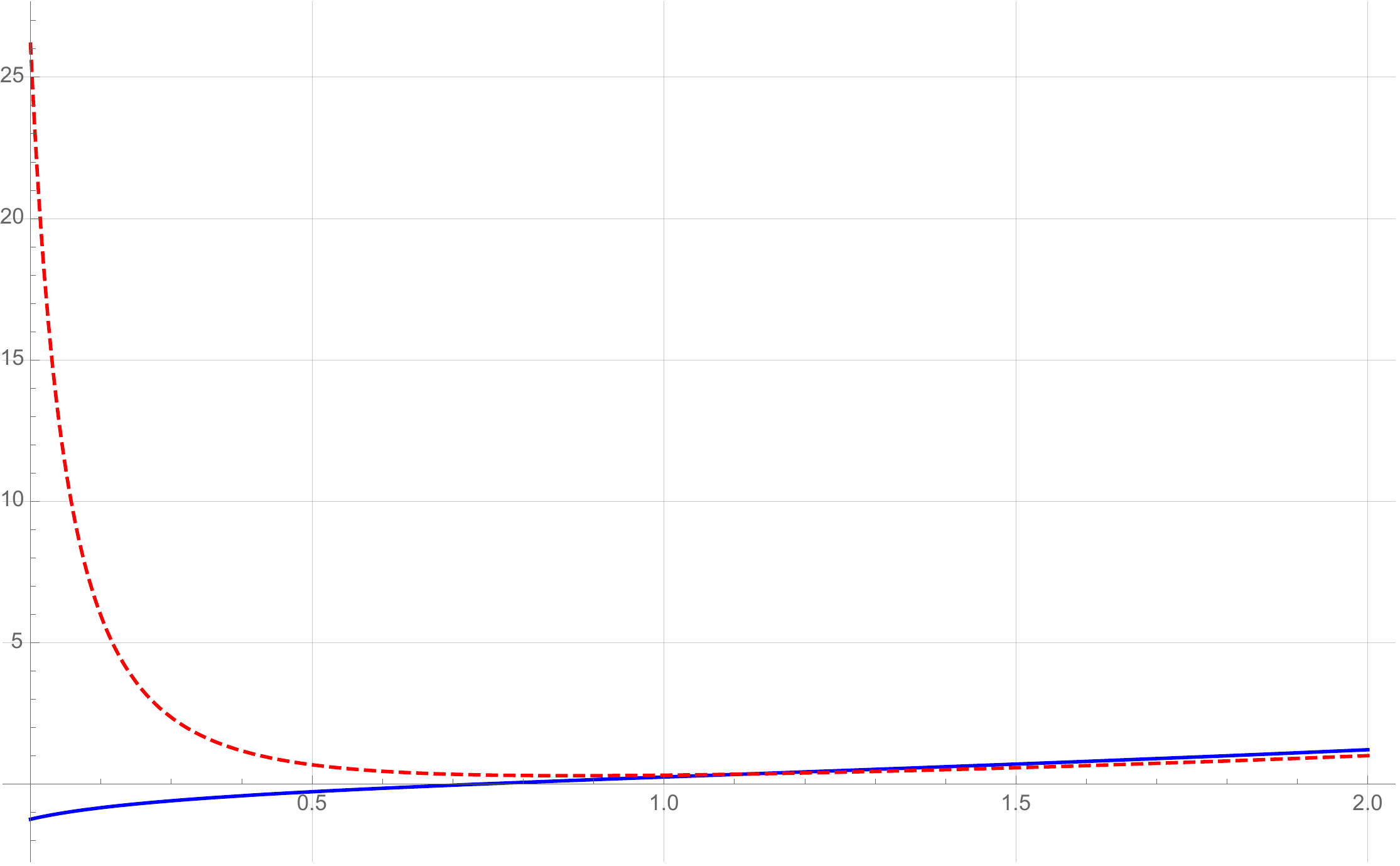}
        \caption{$\beta = 1$ and $N = 6$\,: plot of $\log z_\beta(\sigma)/N^2$}
    \end{subfigure}
    \quad
    \begin{subfigure}[b]{0.475\textwidth}  
        \centering 
        \includegraphics[width=\textwidth]{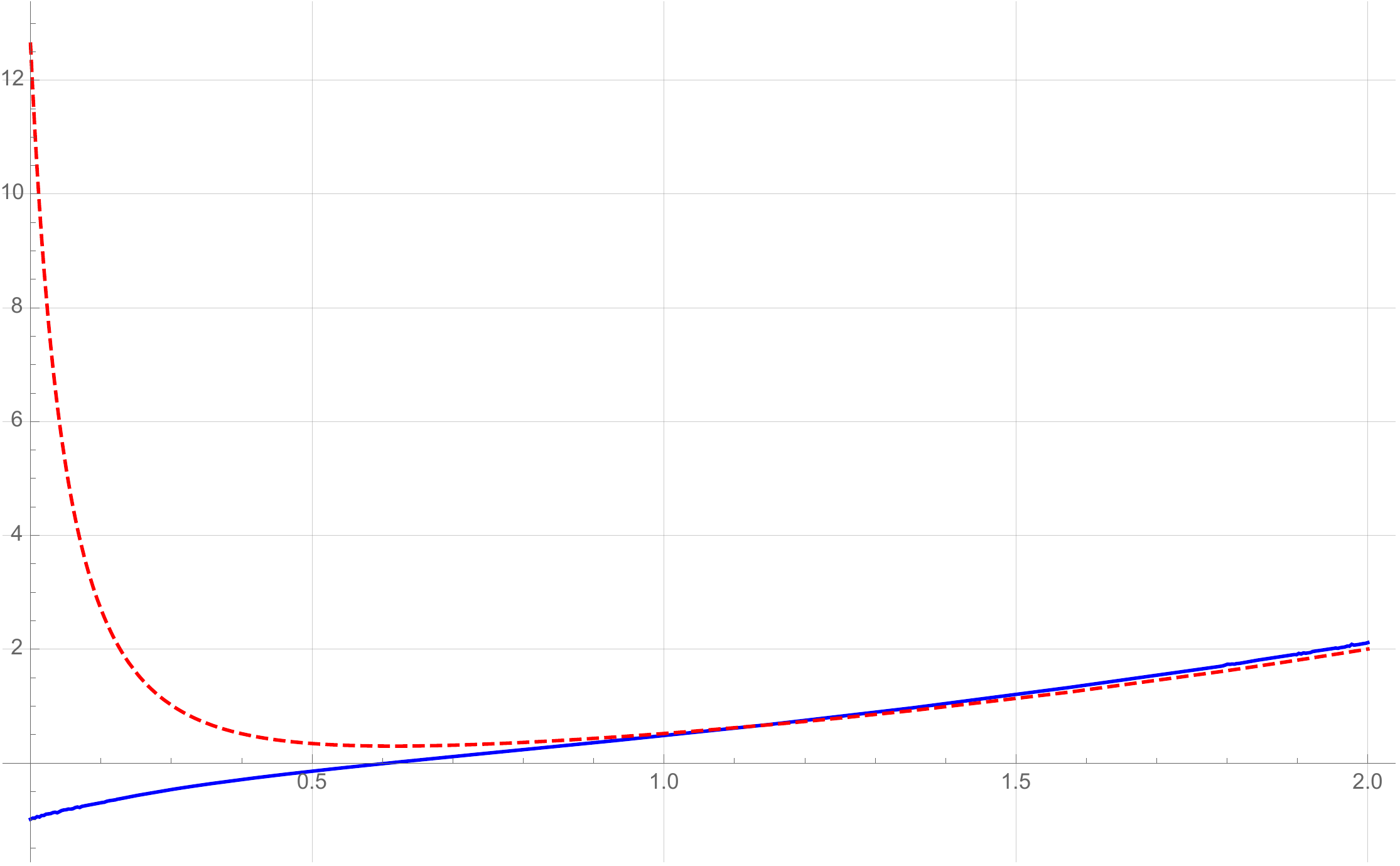}
        \caption{$\beta = 1$ and $N = 12$\,: plot of $\log z_\beta(\sigma)/N^2$}
    \end{subfigure}
    \caption{Proposition \ref{prop:trilog} ($\beta = 1$)\,: (\ref{eq:trilog}) in dashed red and (\ref{eq:pfaff}) in solid blue}
\label{fig:2}
\end{figure*}

It is not equally straightforward to numerically evaluate (\ref{eq:pfaff}). Even at $N=10$, we begin encountering overflow problems for moderate values of $\sigma$ when performing the computations in MATLAB.
Still, as long as these overflow problems do not appear, it is possible to observe a close agreement between (\ref{eq:trilog}) and (\ref{eq:pfaff}), as in the $\beta = 2$ case. This is shown in Figure \ref{fig:2} for $N = 6$ and $12$.

Based on the numerical results summarised in Figures \ref{fig:1} and \ref{fig:2}, it seems possible to use the right-hand side of (\ref{eq:trilog}), multiplied by $N^2$, as a substitute for $\log z_\beta(\sigma)$. While this is only an approximation, it is a highly efficient one, and has the advantage of being given by a rather simple analytic expression. In the $\beta = 1$ case, direct numerical evaluation of $\log z_\beta(\sigma)$ is unstable for larger values of $\sigma$, and (\ref{eq:trilog}) offers a practical way out of this problem. 

\begin{figure*}
    \centering
    \begin{subfigure}[b]{0.475\textwidth}
        \centering
        \includegraphics[width=\textwidth]{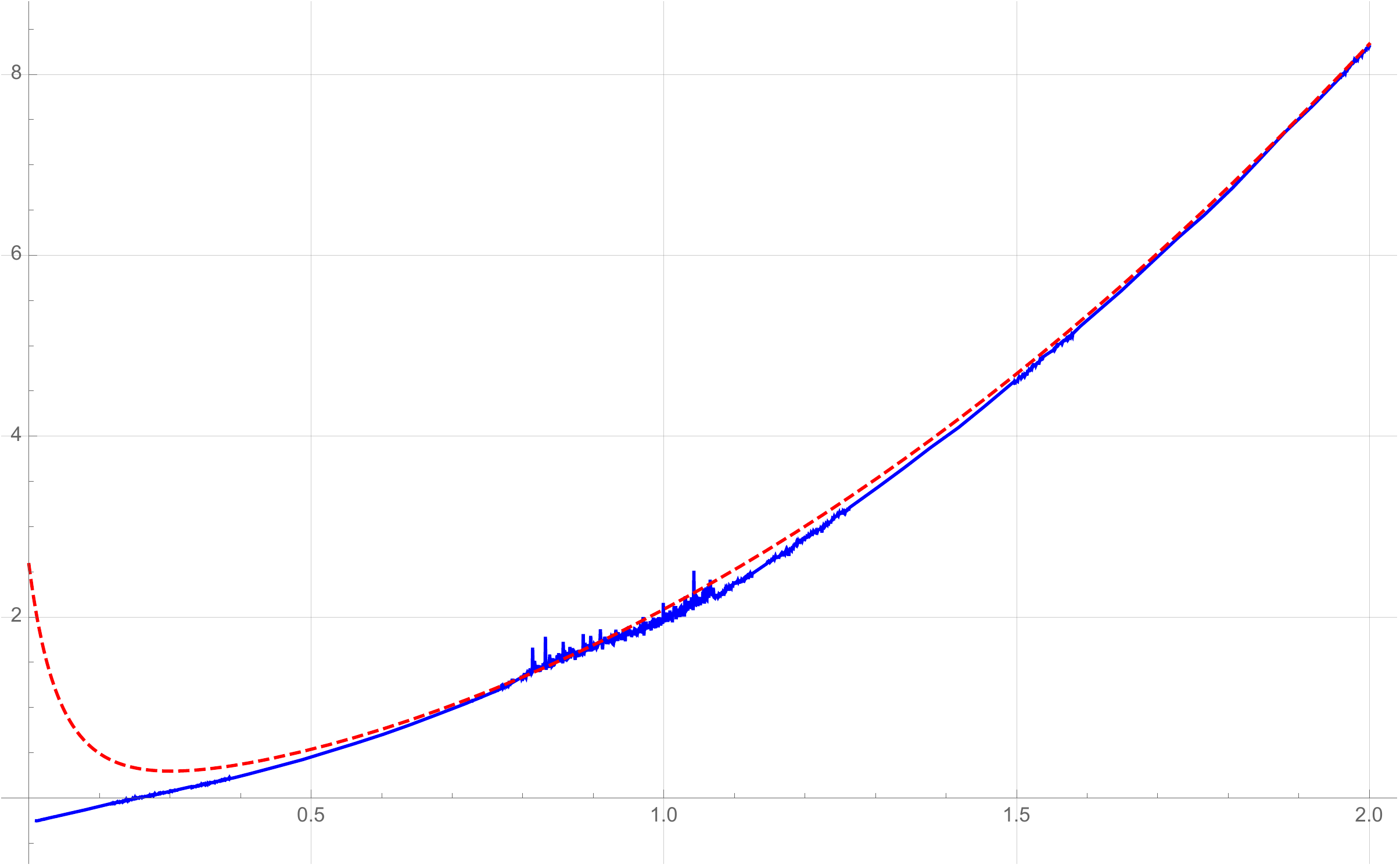}
        \caption{$\beta = 1$ and $N=50$\,: plot of $\log z_\beta(\sigma)/N^2$ up to $\sigma = 2$}
        \label{3a}
    \end{subfigure}
    \quad
    \begin{subfigure}[b]{0.475\textwidth}  
        \centering 
        \includegraphics[width=\textwidth]{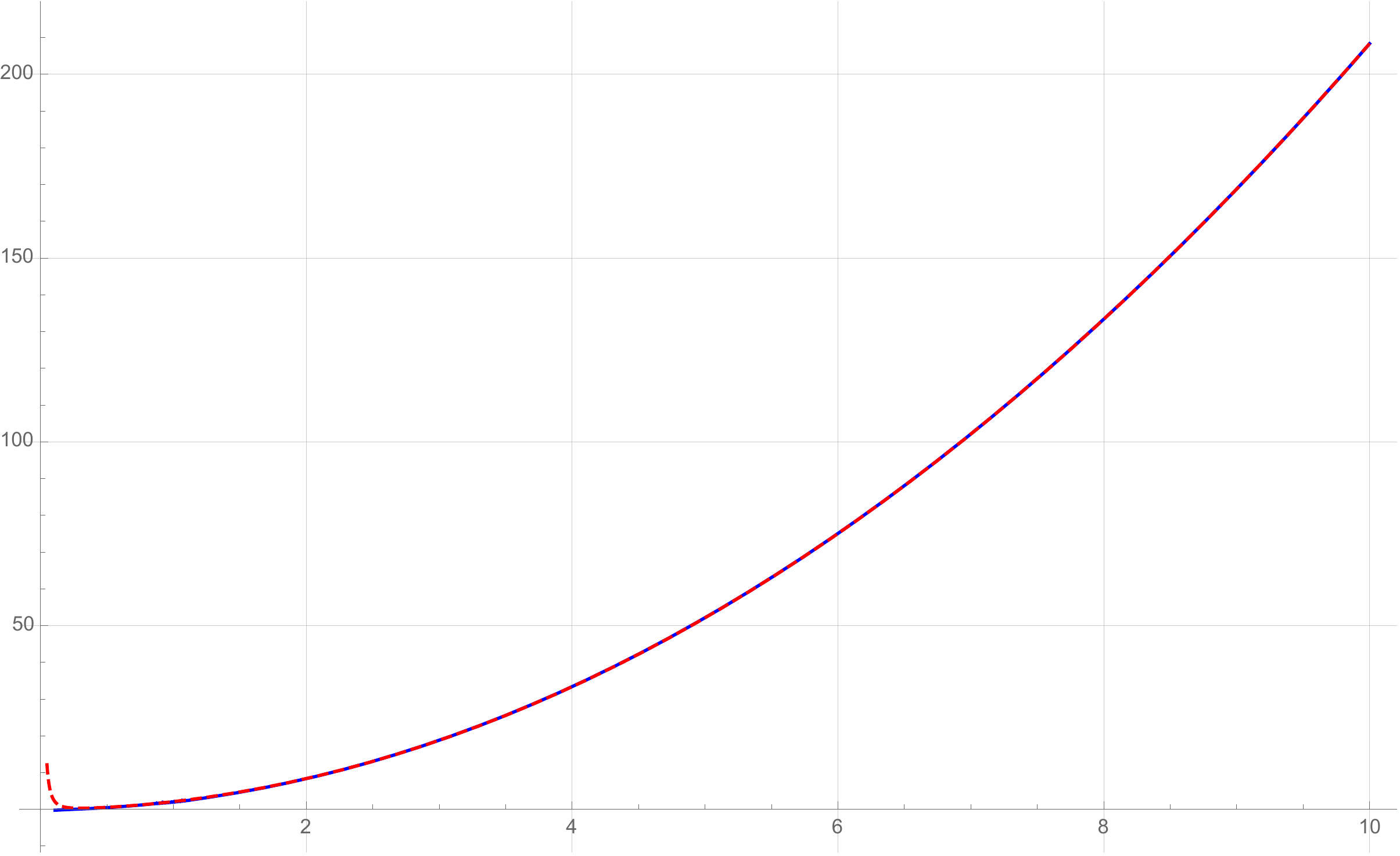}
        \caption{$\beta = 1$ and $N=50$\,: plot of $\log z_\beta(\sigma)/N^2$ up to $
    \sigma=10$}
    \label{3b}
    \end{subfigure}
    \caption{Proposition \ref{prop:trilog} ($\beta = 1$)\,: (\ref{eq:trilog}) in dashed red and (\ref{eq:pfaff}) in solid blue for $N=50$}
\label{fig:3}
\end{figure*}

This can be verified by using symbolic computation (Mathematica), which allows us to extend the evaluation of (\ref{eq:pfaff}) to larger values of $\sigma$ for larger $N$. It is nonetheless associated with the drawback that the resulting curves tend to be artificially non-smooth as in Figure \ref{3a}, due to numerical artifacts. These non-smooth features also exist in the curve depicted in Figure \ref{3b}, but do not appear visible at the given resolution. This behaviour is particularly problematic in the context of Equation (\ref{eq:siglikelihood}), due to the presence of the derivative of $\log z_{\beta}(\sigma)$. Furthermore, symbolic computation of (\ref{eq:pfaff}) becomes exceedingly slow for sufficiently large $N$, and this can only be overcome by relying on approximations such as (\ref{eq:trilog}).

The second set of experiments directly addressed Equation (\ref{eq:siglikelihood}), in the $\beta = 1$ case, for a dimension $N$ ranging between $10$ and $25$. This equation was solved using the Newton method, after its right-hand side was approximated according to (\ref{eq:trilog}). Here, the solution obtained in this way will be denoted $\tilde{\sigma}_{\scriptscriptstyle M}$. This $\tilde{\sigma}_{\scriptscriptstyle M}$ is an approximation of the
maximum-likelihood estimate $\hat{\sigma}_{\scriptscriptstyle M}$ (the exact solution of (\ref{eq:siglikelihood})). If this approximation is any good, $\tilde{\sigma}_{\scriptscriptstyle M}$ should approach the true value of $\sigma$ for sufficiently large $M$ (the number of data points $Y_{\scriptscriptstyle m}$). This can already be observed at $M = 10^3$ when $N = 10$, as reported in the following Table \ref{tab:1}. Each entry in this table gives the average value and standard deviation of $\tilde{\sigma}_{\scriptscriptstyle M}\hspace{0.03cm}$, calculated over $20$ independent trials (average $\pm$ standard deviation). It is clear that increasing $M$ from $10^3$ to $10^4$ only reduces the standard deviation of $\tilde{\sigma}_{\scriptscriptstyle M}\hspace{0.03cm}$, without really affecting its average. 
 \begin{table*}[h]
 \centering
    \begin{tabular}{lccccccc}
    \hline
    true $\sigma$ & 1 & 2 & 3 & 4 & 5 & 6 & 7 \\[0.2cm]
    
   $M = 10^3$     & $ 1.09 \pm 0.00$ & $2.01 \pm 0.01$ & $3.11 \pm 0.30$ & $ 3.56 \pm 0.60$ & $4.58 \pm 0.60$ & $5.13 \pm 0.50$ & $5.60 \pm 0.75$ \\[0.2cm]
   $M = 10^4$     & $1.09 \pm 0.00$ & $2.01 \pm 0.00$ & $3.11 \pm 0.15$ & $3.80 \pm 0.16$ & $4.50 \pm 0.23$ & $5.37  \pm 0.31$ & $5.60 \pm 0.37$ \\[0.1cm]
   \hline
 \end{tabular} 
    \caption{The solution $\tilde{\sigma}_{\scriptscriptstyle M}$ of (\ref{eq:siglikelihood}), for $\beta = 1$ and $N = 10$ (r.h.s. approximated using (\ref{eq:trilog}))}
    \label{tab:1}   
\end{table*}   

When $N = 10$, the right-hand side of (\ref{eq:siglikelihood}) can still be approximated using the Monte Carlo technique mentioned in~\cite{Sa16}\cite{Sa17}. The solution obtained with this approximation will be denoted $\left[\tilde{\sigma}_{\scriptscriptstyle M}\right]_{\mathrm{MC}}\hspace{0.03cm}$. The following Table \ref{tab:2} shows that the estimation error from $\tilde{\sigma}_{\scriptscriptstyle M}$ is quite improved, in comparison with the estimation error from $\left[\tilde{\sigma}_{\scriptscriptstyle M}\right]_{\mathrm{MC}}\hspace{0.03cm}$. Moreover, $\left[\tilde{\sigma}_{\scriptscriptstyle M}\right]_{\mathrm{MC}}$ seems to systematically overestimate the true value of $\sigma$.
 \begin{table*}[h]
 \centering
    \begin{tabular}{lccccccc}
    \hline
    true $\sigma$ & 1 & 2 & 3 & 4 & 5 & 6 & 7 \\[0.2cm]
    
   $M = 10^3$     & $ 1.16 \pm 0.01$ & $2.90 \pm 0.01$ & $5.63 \pm 0.93$ & $ 7.38 \pm 1.35$ & $ 8.35 \pm 1.50 $ & $9.40 \pm 1.11$ & $9.95 \pm 1.12$ \\[0.2cm]
   $M = 10^4$     & $1.16 \pm 0.01$ & $2.90 \pm 0.00$ & $5.40 \pm 0.37$ & $7.00 \pm 0.60$ & $ 9.00 \pm 0.50$ & $9.50  \pm 0.31$ & $9.77 \pm 0.35$ \\[0.1cm]
 \hline
\end{tabular} 
    \caption{The solution $\left[\tilde{\sigma}_{\scriptscriptstyle M}\right]_{\mathrm{MC}}$ of (\ref{eq:siglikelihood}), for $\beta = 1$ and $N = 10$ (MC approximation of r.h.s.)}
    \label{tab:2}   
\end{table*}   

When $N$ is above $20$, Monte Carlo approximation of the right-hand side of (\ref{eq:siglikelihood}) is not feasible anymore, and one is left only with the approximation using (\ref{eq:trilog}). The solution $\tilde{\sigma}_{\scriptscriptstyle M}$ obtained from this approximation is shown in Table \ref{tab:3}, for $N = 20$ and $M = 10^4$. Here (contrary to Table \ref{tab:2}), $\tilde{\sigma}_{\scriptscriptstyle M}$ systematically underestimates the true value of $\sigma$. In fact, identical behavior was observed for $N$ between $20$ and $25$, along with similar values of $\tilde{\sigma}_{\scriptscriptstyle M}\hspace{0.03cm}$. 

A practical means of overcoming this issue would be to include a penalty term into Equation (\ref{eq:siglikelihood}), in order to enforce greater values of its solution. Then, the approximation using (\ref{eq:trilog}) can be successfully implemented for larger dimension $N$, where Monte Carlo approximation becomes useless (the present investigation stopped at $N = 25$, because the sampling algorithms used to generate the data points $Y_{\scriptscriptstyle m}$ could not be taken any further). 

 \begin{table*}[h!]
 \centering
    \begin{tabular}{lccccccc}
    \hline
    true $\sigma$ & 1 & 2 & 3 & 4 & 5 & 6 & 7 \\[0.2cm]
    
   $M = 10^4$     & $ 1.02 \pm 0.00$ & $1.51 \pm 0.00$ & $1.87 \pm 0.01$ & $ 2.35 \pm 0.02$ & $ 2.82 \pm 0.02$ & $3.24 \pm 0.03$ & $3.63 \pm 0.02$ \\[0.1cm]
   \hline
 \end{tabular} 
    \caption{The solution $\tilde{\sigma}_{\scriptscriptstyle M}$ of (\ref{eq:siglikelihood}), for $\beta = 1$ and $N = 20$ (r.h.s. approximated using (\ref{eq:trilog}))}
    \label{tab:3}   
\end{table*} 

In conclusion, the present section has demonstrated the new approximation of $\log z_\beta(\sigma)$, based on (\ref{eq:trilog}), significantly improves on the existing Monte Carlo approximation, and also extends it to higher dimensions. Moreover, (\ref{eq:trilog}) makes up for the numerical instability and computational cost of exact formulae such as (\ref{eq:pfaff}), and its straightforward analytic form makes it possible to employ the Newton method in solving (\ref{eq:siglikelihood}). This is preferable to employing a grid search based on (\ref{eq:pfaff}), whose performance is restricted by grid resolution and size.  

In upcoming work, the experiments conducted in this section will be generalised to more realistic learning models which may be applied to real-world data (for example, the mixture models considered in~\cite{zanini2018}).%

\section{Block-Toeplitz covariance matrices} \label{sec:siegel}
Similar to the spaces $\mathcal{P}^{\hspace{0.02cm}\scriptscriptstyle \beta}_{\scriptscriptstyle N}$ introduced in Section \ref{sec:background}, consider the spaces $\mathcal{D}^{\hspace{0.02cm}\scriptscriptstyle \beta}_{\scriptscriptstyle N\hspace{0.03cm}}$, defined as follows.
\begin{itemize}
    \item $\mathcal{D}^{\hspace{0.02cm}\scriptscriptstyle 2}_{\scriptscriptstyle N}$ is the space of $N \times N$ matrices with complex entries, whose operator norm is $< 1$.
    \item $\mathcal{D}^{\hspace{0.02cm}\scriptscriptstyle 1}_{\scriptscriptstyle N}$ is the space of $N \times N$ symmetric matrices with complex entries, whose operator norm is $< 1$.
\end{itemize}
Here, operator norm means the largest singular value. The space $\mathcal{D}^{\hspace{0.02cm}\scriptscriptstyle 2}_{\scriptscriptstyle N}$
will be called the Hermitian Siegel domain, and $\mathcal{D}^{\hspace{0.02cm}\scriptscriptstyle 1}_{\scriptscriptstyle N}$ the symmetric Siegel domain. The focus will be restricted to these $\beta = 1$ or $2$ cases, since they are closely related to time-series analysis and signal processing (the $\beta = 4$ case is similar, but involves quaternion matrices). 

Specifically, assume a wide-sense stationary $N$-variate time series of length $T$ is described by its autocovariance matrix $\Gamma$. Wide-sense stationarity implies $\Gamma$ has a block-Toeplitz structure, with $T \times T$ blocks of size $N \times N$. When solving an optimal prediction problem, one may apply a multidimensional Szegö-Levinson algorithm to the autocovariance $\Gamma$~\cite{jeuris1}\cite{yannthesis}, and obtain a family of matrices $(\Gamma_{\scriptscriptstyle 0},\Omega_{\scriptscriptstyle 1},\ldots, \Omega_{\scriptscriptstyle T-1})$, where $\Gamma_{\scriptscriptstyle 0}$ is the zero-lag autocovariance of the original time series (this is a complex covariance matrix, $\Gamma_{\scriptscriptstyle 0} \in \mathcal{P}^{\hspace{0.02cm}\scriptscriptstyle 2}_{\scriptscriptstyle N\hspace{0.03cm}}$), and each $\Omega_t$ is a so-called matrix reflection coefficient, which belongs to $\mathcal{D}^{\hspace{0.02cm}\scriptscriptstyle 2}_{\scriptscriptstyle N\hspace{0.03cm}}$. If the autocovariance $\Gamma$ has Toeplitz blocks, then each $\Omega_t$ moreover belongs to $\mathcal{D}^{\hspace{0.02cm}\scriptscriptstyle 1}_{\scriptscriptstyle N}$ (which is a subspace of  $\mathcal{D}^{\hspace{0.02cm}\scriptscriptstyle 2}_{\scriptscriptstyle N}$).   

Riemannian geometry enters this picture in the following way~\cite{jeuris1}\cite{yannthesis}. Consider an information metric  
on the space of block-Toeplitz autocovariance matrices $\Gamma$, equal to the Hessian of the entropy function $S(\Gamma) = \log\det(\Gamma)$. In terms of the new coordinates $(\Gamma_{\scriptscriptstyle 0},\Omega_{\scriptscriptstyle 1},\ldots, \Omega_{\scriptscriptstyle T-1})$, this information metric is a direct product of the affine-invariant metric (\ref{eq:aff-inv}) on the first coordinate $\Gamma_{\scriptscriptstyle 0\hspace{0.03cm}}$, and of a scaled copy of the Siegel domain metric on each of the remaining coordinates $\Omega_t\hspace{0.03cm}$ (see formula (3.11) in~\cite{jeuris1}). Dropping the subscript $t$, this Siegel domain metric is given by,
\begin{equation}\label{eq:siegelmetric}   
\!\!\langle u\hspace{0.02cm},v\rangle_{\scriptscriptstyle \Omega} =  \Re\,\mathrm{tr}\left[ \left(\mathrm{I}_{\scriptscriptstyle N} - \Omega\Omega^\dagger\right)^{\scriptscriptstyle -1}u\left(\mathrm{I}_{\scriptscriptstyle N} - \Omega\Omega^\dagger\right)^{\scriptscriptstyle -1}v^\dagger\right]
\end{equation}
for $u\hspace{0.02cm},v \in T_{\scriptscriptstyle \Omega}\mathcal{D}^{\hspace{0.02cm}\scriptscriptstyle \beta}_{\scriptscriptstyle N}\hspace{0.03cm}$,
where the tangent space $T_{\scriptscriptstyle \Omega}\mathcal{D}^{\hspace{0.02cm}\scriptscriptstyle \beta}_{\scriptscriptstyle N}$ is isomorphic to the space of $N \times N$ complex matrices if $\beta = 2$, and to the space of $N \times N$ symmetric complex matrices if $\beta = 1$ (this isomorphism shows that the dimension of $\mathcal{D}^{\hspace{0.02cm}\scriptscriptstyle \beta}_{\scriptscriptstyle N}$ equals $2NN_\beta$). The geodesic distance associated with the Siegel metric (\ref{eq:siegelmetric}) has the following expression~\cite{jeuris1}\cite{yannthesis},
\begin{equation} \label{eq:siegeldistance}
 d^{\hspace{0.03cm} 2}(\Xi\hspace{0.02cm},\Omega) = \mathrm{tr}\!\left(\mathrm{arctanh}^2\!\left(R^{\scriptscriptstyle \frac{1}{2}}(\Xi\hspace{0.02cm},\Omega)\right)\right)
 \hspace{0.4cm}
\Xi\hspace{0.02cm},\Omega \in \mathcal{D}^{\hspace{0.02cm}\scriptscriptstyle \beta}_{\scriptscriptstyle N}
\end{equation} 
where 
\begin{equation*}
R(\Xi\hspace{0.02cm},\Omega) = (\Xi - \Omega)(\mathrm{I}_{\scriptscriptstyle N} - \Omega^\dagger\Xi)^{\scriptscriptstyle -1}(\Xi^\dagger - \Omega^\dagger)(\mathrm{I}_{\scriptscriptstyle N} - \Omega\hspace{0.02cm}\Xi^\dagger)^{\scriptscriptstyle -1}
\end{equation*}
is called the matrix cross-ratio.

Riemannian Gaussian distributions on the space $\mathcal{D}^{\hspace{0.02cm}\scriptscriptstyle \beta}_{\scriptscriptstyle N}$ are given by their probability density function $p(\Omega|\bar{\Omega},\sigma)$ which is of the same form as in (\ref{eq:gdensity}), but with the distance $d(\Omega\hspace{0.02cm},\bar{\Omega})$ determined by (\ref{eq:siegeldistance}), and the normalising factor
\begin{equation} \label{eq:Zsiegel}
Z(\sigma) = \int_{\mathcal{D}^{\hspace{0.02cm}\scriptscriptstyle \beta}_{\scriptscriptstyle N}}\,\exp\left[-\frac{d^{\hspace{0.03cm} 2}(\Omega\hspace{0.02cm},\bar{\Omega})}{2\sigma^2}\right]dv(\Omega)
\end{equation}
with respect to the Riemannian volume $dv(\Omega) = \det(\mathrm{I}_{\scriptscriptstyle N} - \Omega\Omega^\dagger)^{-2N_\beta}\lbrace d\Omega\rbrace$. Here, $\lbrace d\Omega \rbrace =\prod_{ij} \Re\,d\Omega_{ij}\Im\,d\Omega_{ij}$ where $\Re$ and $\Im$ denote the real and imaginary parts (the product is over $i \leq j$ if $\beta = 1$ and over all $i\hspace{0.02cm},j$ if $\beta = 2$).

As shown in~\cite{Sa17}, a Riemannian Gaussian distribution on the space of block-Toeplitz covariance matrices $\Gamma$ is just a product of independent Riemannian Gaussian distributions, one for each coordinate $(\Gamma_{\scriptscriptstyle 0},\Omega_{\scriptscriptstyle 1},\ldots, \Omega_{\scriptscriptstyle T-1})$. For $\Gamma_{\scriptscriptstyle 0}\hspace{0.03cm}$, this is a Riemannian Gaussian distribution on $\mathcal{P}^{\hspace{0.02cm}\scriptscriptstyle 2}_{\scriptscriptstyle N\hspace{0.03cm}}$, already considered in Section \ref{sec:background}. Thus, to understand Riemannian Gaussian distributions of block-Toeplitz covariance matrices, it only remains to study Riemannian Gaussian distributions on the Siegel domain $\mathcal{D}^{\hspace{0.02cm}\scriptscriptstyle \beta}_{\scriptscriptstyle N}$. 

In Section \ref{sec:contribution}, it was seen that Riemannian Gaussian distributions on $\mathcal{P}^{\hspace{0.02cm}\scriptscriptstyle \beta}_{\scriptscriptstyle N}$ are equivalent to log-normal matrix ensembles. On the other hand, Riemannian Gaussian distributions on $\mathcal{D}^{\hspace{0.02cm}\scriptscriptstyle \beta}_{\scriptscriptstyle N}$ are equivalent to ``acosh-normal" ensembles, which will be described in Proposition \ref{prop:acosh} below. 

Note first that each matrix $\Omega \in \mathcal{D}^{\hspace{0.02cm}\scriptscriptstyle \beta}_{\scriptscriptstyle N}$ 
can be factorised in the following way
\begin{equation} \label{eq:takagi}
 \Omega = U\lambda V^\dagger \;\; (\text{if } \beta = 2) \hspace{0.5cm}
 \Omega = U\lambda\hspace{0.02cm}U^T \;\;(\text{if } \beta = 1)
\end{equation}
where $U$ and $V$ are unitary, $^T$ denotes the transpose, and $\lambda$ is diagonal, with diagonal elements $\lambda_ i =  \mathrm{tanh}(r_i)$ for $(r_{\scriptscriptstyle 1},\ldots,r_{\scriptscriptstyle N}) \in \mathbb{R}^N$. The $\beta = 2$ case follows from the singular value decomposition of $\Omega$, and the $\beta = 1$ case follows from the Takagi decomposition of $\Omega$~\cite{edelman}. In either case, $\lambda_ i$ is of the form $\mathrm{tanh}(r_i)$ because the singular values of $\Omega$ are all $< 1$.
\begin{prop} \label{prop:acosh}
  Let $\Omega$ follow a Riemannian Gaussian distribution on the Siegel domain $\mathcal{D}^{\hspace{0.02cm}\scriptscriptstyle \beta}_{\scriptscriptstyle N}$ with $\bar{\Omega} = \mathrm{0}_{\scriptscriptstyle N}$ ($N \times N$ zero matrix). If $\Omega$ is factorised as in (\ref{eq:takagi}), then $U$ and $V$ are uniformly distributed on the unitary group $U(N)$. Moreover, if $x_i = \cosh(2r_i)$ then $(x_{\scriptscriptstyle 1},\ldots,x_{\scriptscriptstyle N})$ have joint probability density function 
 \begin{equation} \label{eq:acoshnormal2}
  p(x|\sigma) \propto |V(x)|^\beta\prod_i \rho(x_i\hspace{0.02cm},8\sigma^2)  
\end{equation}
where $(x_{\scriptscriptstyle 1},\ldots,x_{\scriptscriptstyle N}) \in (1\hspace{0.02cm},\infty)^N$, $V(x) = \prod_{i<j}(x_j-x_i)$, and $\rho(x\hspace{0.02cm},k) = \exp(-\mathrm{acosh}^2(x)/k)$ is the ``acosh-normal" weight function.
\end{prop}
Note that the assumption that the centre-of-mass parameter $\bar{\Omega}$ is equal to $\mathrm{0}_{\scriptscriptstyle N}$ does not entail any loss of generality. The transformation~\cite{yannthesis}
\begin{align} \label{eq:moebius1}
\nonumber  \Omega \longmapsto \Psi(\Omega) = \\  (\mathrm{I}_{\scriptscriptstyle N} - \bar{\Omega}\bar{\Omega}^\dagger)^{\scriptscriptstyle -\frac{1}{2}}(\Omega - \bar{\Omega})(\mathrm{I}_{\scriptscriptstyle N} - \bar{\Omega}^\dagger\Omega)^{\scriptscriptstyle -1}
  (\mathrm{I}_{\scriptscriptstyle N} - \bar{\Omega}^\dagger\bar{\Omega})^{\scriptscriptstyle \frac{1}{2}}
\end{align}
maps $\bar{\Omega}$ to $\mathrm{0}_{\scriptscriptstyle N}$, while preserving the Siegel domain metric (\ref{eq:siegelmetric}) and the associated distance and Riemannian volume. Thus, if $\Omega$ follows a Riemannian Gaussian density $p(\Omega|\bar{\Omega},\sigma)$, it is enough to replace $\Omega$ by $\Psi(\Omega)$, which will have a Riemannian Gaussian density with $\bar{\Omega} = \mathrm{0}_{\scriptscriptstyle N}$ and with the same $\sigma$.

Proposition \ref{prop:acosh} implies that the transformed singular values $x_i = \cosh(2r_i)$ follow the classical eigenvalue distribution of an orthogonal ($\beta = 1$) or unitary ($\beta = 2$) matrix ensemble with ``acosh-normal" weight function. Therefore, in particular, the normalising factor $Z(\sigma)$ of  (\ref{eq:Zsiegel}) reduces to a multiple integral (compare to (\ref{eq:Z2}) and (\ref{eq:z2}))
\begin{equation} \label{eq:zsiegel}
z_\beta(\sigma) = \frac{1}{N!}\int_{(1\hspace{0.02cm},\infty)^N}\,|V(x)|^\beta\,\omega(dx)
\end{equation}
where $\omega(dx) = \prod\, \rho(dx_i)$, with $\rho(dx_i) = \rho(x_i\hspace{0.02cm},8\sigma^2)\hspace{0.02cm}dx_i$ for the weight function $\rho(x\hspace{0.02cm},k) = \exp(-\mathrm{acosh}^2(x)/k)$.

Thus, Riemannian Gaussian distributions on the Siegel domain $\mathcal{D}^{\hspace{0.02cm}\scriptscriptstyle \beta}_{\scriptscriptstyle N}$ can be studied by applying the tools of random matrix theory to new ``acosh-normal" matrix ensembles. Hopefully, in the near future, this will lead to similar results to the ones obtained with log-normal ensembles in Section \ref{sec:contribution}, paving the way to learning from datasets of high-dimensional block-Toeplitz covariance matrices.

\appendices

\section{Proofs of Propositions \ref{prop:log-normal} to \ref{prop:acosh}} \label{app:proofs}

\subsection{Proof of Propositon \ref{prop:log-normal}}
\noindent 

If $Y$ follows the Riemannian Gaussian density (\ref{eq:gdensity}) with $\bar{Y} = \mathrm{I}_{\scriptscriptstyle N}\hspace{0.03cm}$, then 
\begin{equation} \label{eq:prooflnr1}
\mathbb{P}(Y \in B) = \left(Z(\sigma)\right)^{-1}\,\int_B\,
\exp\left[-\frac{d^{\hspace{0.03cm} 2}(Y\hspace{0.02cm},\mathrm{I}_{\scriptscriptstyle N})}{2\sigma^2}\right]dv(Y)
\end{equation}
Let $(y_{\scriptscriptstyle 1},\ldots,y_{\scriptscriptstyle N})$ denote the eigenvalues of $Y$. Using (\ref{eq:aff-dist}) and the fact that $dv(Y) = \det(Y)^{-N_\beta}\lbrace dY\rbrace$, (\ref{eq:prooflnr1}) becomes
\begin{align}\label{eq:prooflnr2}
\nonumber \mathbb{P}(Y \in B) = \\ \left(Z(\sigma)\right)^{-1}\,\int_B\,\left(
\prod^N_{i=1}\exp\left[-\frac{\log^2(y_i)}
{2\sigma^2}\right]y^{-N_\beta}_i\right)\lbrace dY\rbrace
\end{align}
Recall that $X = e^{N_\beta\sigma^2}\hspace{0.02cm}Y$. Accordingly,
if $(x_{\scriptscriptstyle 1},\ldots,x_{\scriptscriptstyle N})$ are the eigenvalues of $X$, an elementary calculation yields
$$
\begin{array}{r}
\exp\!\left[-\frac{\log^2(y_i)}
{2\sigma^2}\right]y^{-N_\beta}_i \!=\! \exp[N^2_\beta(\sigma^2/2)]
\exp\!\left[-\frac{\log^2(x_i)}
{2\sigma^2}\right]
\end{array}
$$
Therefore, (\ref{eq:prooflnr2}) can be written
\begin{align}\label{eq:prooflnr3}
\nonumber \mathbb{P}(Y \in B) = \exp[NN^2_\beta(\sigma^2/2)]\times\left(Z(\sigma)\right)^{-1}\,\\ \int_B\,\mathrm{etr}\left[-\frac{\log^2(X)}{2\sigma^2}\right] \lbrace dY\rbrace
\end{align}
To conclude, it is enough to use once more the definition $X = e^{N_\beta\sigma^2}\hspace{0.02cm}Y$. This implies
$$
\lbrace dY \rbrace = \exp[-2NN^2_\beta(\sigma^2/2)]\hspace{0.03cm}\lbrace dX\rbrace
$$
Thus, using the fact that
$$
\mathbb{P}(X \in B) =
\mathbb{P}\left(e^{N_\beta\sigma^2}\hspace{0.02cm}Y \in B\right)
$$
and changing the variable of integration from $Y$ to $X$ in (\ref{eq:prooflnr3}), it follows that
\begin{align}
\nonumber \mathbb{P}(X \in B) =
\exp[-NN^2_\beta(\sigma^2/2)]\times\left(Z(\sigma)\right)^{-1}\,\\ \nonumber \int_B\,\mathrm{etr}\left[-\frac{\log^2(X)}{2\sigma^2}\right] \lbrace dX\rbrace
\end{align}
as required in (\ref{eq:lognormal1}).

\subsection{Proof of Proposition \ref{prop:trilog}}
The $\beta = 2$ case follows from (\ref{eq:swp}), by an elementary calculation, after noting that
$$
    \frac{1}{N^2}\sum^{N-1}_{n=1}(N-n)\log\!\left(1 - e^{-n\hspace{0.02cm}\sigma^2}\right)
$$
is a Riemann sum for the improper integral
$$
\int_0^1(1-x)\log\!\left(1-e^{-tx}\right)\hspace{0.02cm}dx = \frac{\mathrm{Li}_{\scriptscriptstyle 3}(1) -\mathrm{Li}_{\scriptscriptstyle 3}(e^{-t})}{t^2}
$$
For other values of $\beta$, the result will be obtained from the reasoning presented in the proof of Proposition \ref{prop:empirical}, based on the scaling equation (\ref{scaling result}).

\subsection{Proof of Proposition \ref{prop:empirical}}
The $\beta = 2$ case has already been proved in~\cite{habilitation}. To deal with the general case, write (\ref{eq:z2}) under the form $z_{\beta}(\sigma)=c_{\beta}(\sigma)\,I_{\beta}(\sigma)$, where
$$
c_{\beta}(\sigma)=\frac{1}{N!}\exp\left[-NN^2_\beta(\sigma^2/2)\right]
$$
and $I_\beta(\sigma)$ is the multiple integral
\begin{align}
\nonumber I_\beta(\sigma) = \int_{\mathbb{R}^N_+}\exp\left[-\frac{1}{2\sigma^2}\sum^N_{i=1} \log^2(x_i)\right. \,+ \\
\nonumber \left. \beta\,\sum_{i<j} \log|x_j - x_i|\right]dx 
\end{align}
Now, if $\mu_{\scriptscriptstyle N} = (1/N)\sum^N_{i=1}\,\delta_{x_i}$ is the empirical distribution of $x_i$ (where $\delta_{x_i}$ denotes the Dirac measure at $x_i$), then 
\begin{equation} \label{eq:ibeta22}
I_\beta(\sigma) = \int_{\mathbb{R}^N_+}\exp\left[-N^2E_{\beta}(\mu_{\scriptscriptstyle N},t)\right]dx 
\end{equation}
where $E_{\beta}(\mu,t)$ is the so-called energy functional
\begin{align} \label{eq:energy1}
\nonumber E_{\beta}(\mu,t) = \frac{1}{2t}\int_{\mathbb{R}_+}\log^2(x)\hspace{0.02cm}\mu(dx) \,- \\\beta \int_{\mathbb{R}_+}\int_{\mathbb{R}_+} \log|y-x|\hspace{0.02cm}\mu(dx)\mu(dy)
\end{align}
defined for any probability distribution $\mu$ on $\mathbb{R}_+\hspace{0.02cm}$.
This energy functional satisfies the following scaling equation
\begin{equation} \label{scaling result}
    E_{\beta}(\mu,t)=(\beta/2)\hspace{0.03cm}E_{2}\!\left(\mu,(\beta/2)\hspace{0.03cm} t\right)
\end{equation}
which is easily obtained after dividing (\ref{eq:energy1}) by $\beta/2$. The proof can now be completed by applying to (\ref{eq:ibeta22}) and (\ref{scaling result}) the arguments in~\cite{deift}. First~\cite{deift} (Corollary 6.90, Page 155), 
$$
\frac{1}{N^2}\hspace{0.02cm}\log I_\beta(\sigma) \longrightarrow - E_{\beta}(\mu^*,t)
$$
where $\mu^* = \mu^*_\beta(t)$ is the so-called equilibrium distribution, the unique minimiser of the energy (\ref{eq:energy1}) among probability distributions on $\mathbb{R}_+\hspace{0.02cm}$. From the scaling equation (\ref{scaling result}), it is now clear
\begin{equation} \label{eq:frenergy}
\frac{1}{N^2}\hspace{0.02cm}\log I_\beta(\sigma) \longrightarrow - (\beta/2)\hspace{0.03cm}E_{2}\!\left(\mu^*,(\beta/2)\hspace{0.03cm} t\right) 
\end{equation}
Therefore, by adding the limit of $\log c_\beta(\sigma)/N^2$, 
\begin{equation} \label{eq:frenergy1}
\frac{1}{N^2}\hspace{0.02cm}\log z_\beta(\sigma) \longrightarrow - (\beta/2)\left((\beta t/4) + \hspace{0.03cm}E_{2}\!\left(\mu^*,(\beta/2)\hspace{0.03cm} t\right)\right)
\end{equation}
However, when $\beta = 2$, it is already known that this limit is $\Phi(t)$, where $\Phi$ is defined in (\ref{eq:trilog}). This provides $E_{2}\!\left(\mu^*, t\right) = -t/2 - \Phi(t)$, which can be replaced back into (\ref{eq:frenergy1}), yielding the general case of Proposition \ref{prop:trilog}. 

To prove Proposition \ref{prop:empirical}, it is possible to use the scaling equation (\ref{scaling result}), once more. From~\cite{deift} (Section 6.4), if $\hat{\mu}_{\scriptscriptstyle N}$ is defined as in (\ref{eq:R1}), but with the $x_i$ instead of the $y_i\hspace{0.02cm}$, then $\hat{\mu}_{\scriptscriptstyle N}$ converges weakly to the equilibrium distribution $\mu^* = \mu^*_\beta(t)$. On the other hand, the scaling equation (\ref{scaling result}) implies $\mu^* = \mu^*_2(\beta t/2)$, and $\mu_2(t)$ is already known to be the image of the distribution with density $n(y|t)$, defined as in (\ref{eq:empirical}), under the change of variables $x = e^{\frac{\beta}{2}t}\hspace{0.02cm}y$. This provides $\mu^*$ for any value of $\beta$, and the proposition then follows by changing the variables back from $x$ to $y$. 

\subsection{Proof of Proposition \ref{prop:acosh}} 
The proof relies on the general theory of Gaussian distributions on Riemannian symmetric spaces, as outlined in~\cite{Sa17}\cite{habilitation}. Here, to keep the proof self-contained, it will be helpful to briefly recall certain aspects of this theory (for a more detailed, in-depth discussion, the reader is referred to~\cite{habilitation}, Sections 1.9 and 3.3).

A Gaussian distribution on a Riemannian symmetric space $M$ is defined by its probability density function
\begin{equation} \label{eq:symspace_gdensity}
 p(x|\bar{x},\sigma) = \left(Z(\sigma)\right)^{-1}\exp\left[-\frac{d^{\hspace{0.03cm} 2}(x\hspace{0.02cm},\bar{x})}{2\sigma^2}\right] \hspace{1cm} x \in M
\end{equation}
with respect to the Riemannian volume $dv(x)$ on $M$, with $\bar{x} \in M$ and $\sigma > 0$. 

It is always assumed $M$ is a Riemannian symmetric space of non-positive curvature, associated to a symmetric pair $(G,K)$. Precisely, $G$ is a connected Lie group which acts transitively and isometrically on $M$ (this action is denoted $x \mapsto g\cdot x$, for $g \in G$ and $x \in M$), and $K$ is a compact subgroup of $G$, made up of those elements $k \in G$ which fix a certain point $o \in M$ (that is, $k \in K$ if and only if $k\cdot o = o$). 

For a concrete understanding of Gaussian distributions on $M$, it is necessary to take a closer look at the Lie algebras of $G$ and $K$, denoted $\mathfrak{g}$ and $\mathfrak{k}$, respectively. These are related together by the so-called Cartan decomposition, $\mathfrak{g} = \mathfrak{k} + \mathfrak{p}$ (direct sum), where the subspace $\mathfrak{p}$ of $\mathfrak{g}$ can be identified with the tangent space $T_oM$ (tangent space to $M$ at at $o$). In terms of this decomposition, the main construction needed for the proof can be described as follows. 

Without any loss of generality, $G$ and $K$ are taken to be matrix Lie groups, and $\mathfrak{g}$ and $\mathfrak{k}$ their matrix Lie algebras. Let $\mathfrak{a}$ be a maximal abelian subspace of $\mathfrak{p}$ (that is, all the matrices $a \in \mathfrak{a}$ commute with each other). Then, each matrix $u \in \mathfrak{p}$ can be written under the form $u = k a k^{-1}$, where $k \in K$ and $a \in \mathfrak{a}$ (this factorisation is the general template for the fifty three matrix factorisations outlined in~\cite{edelman}). Moreover, any $x \in M$ admits a representation
\begin{equation}\label{eq:polarcoordinates}
  x = \exp\left(k a k^{-1}\right)\cdot o \hspace{1cm} k \in K \text{ and } a \in \mathfrak{a}
\end{equation}
where $\exp$ denotes the matrix exponential. Incidentally, this representation is not unique, but for almost all $x \in M$ (for all $x \in M$, except a subset of zero volume) $x$ has exactly $|W|$ couples $(k\hspace{0.02cm},a)$ which satisfy (\ref{eq:polarcoordinates}), where $|W|$ is the number of elements of the Weyl group of the symmetric couple $(G,K)$. 

It is now possible to state the following Lemma \ref{lem:rgdsym}~\cite{Sa17}\cite{habilitation}, which will yield the entire proof, by direct application.
\begin{lem} \label{lem:rgdsym}
Let $x$ follow the Gaussian density (\ref{eq:symspace_gdensity}) on $M$, with $\bar{x} = o\hspace{0.03cm}$. If $x$ is represented as in (\ref{eq:polarcoordinates}), then $k$ is uniformly distributed on the compact group $K$. Moreover, $a$ has the following probability density function on $\mathfrak{a}$ (recall that $\mathfrak{a}$ is a real vector space, so the density is with respect to the usual Lebesgue measure on $\mathfrak{a}$)
\begin{equation} \label{eq:aadensity}
  p(a|\sigma) \propto \exp\left[-\frac{\Vert a\Vert^2_o}{2\sigma^2}\right]\prod_{\rho \in \Delta_+} \left(\sinh\left|\rho(a)\right|\right)^{m_\rho}
\end{equation}
where $\Vert a \Vert_o$ is the Riemannian norm of $a$ (since $a \in \mathfrak{p}$, it can be identified with a vector in $T_oM$), and $\Delta_+$ is a set of positive roots $\rho$ on $\mathfrak{a}$, with respective multiplicities $m_\rho$ (each $\rho \in \Delta_+$ is a certain linear function $\rho:\mathfrak{a} \rightarrow \mathbb{R}$).
\end{lem}
Recall that positive roots $\rho : \mathfrak{a} \rightarrow \mathbb{R}$ are any set of linear functions on $\mathfrak{a}$ such that, for any $a \in \mathfrak{a}$, the eigenvalues of the linear operator $\mathrm{ad}_a : \mathfrak{g} \rightarrow \mathfrak{g}$, given by $\mathrm{ad}_a(\xi) = [a\hspace{0.02cm},\xi]$ (this is $a\hspace{0.02cm}\xi - \xi a$), are equal to $\pm \rho(a)$ with respective multiplicities $m_\rho\hspace{0.03cm}$. 

The proof may now begin in earnest. It merely consists in identifying $G$, $K$, $\mathfrak{g}$, $\mathfrak{k}$, $\mathfrak{p}$, and $\mathfrak{a}$, for each Siegel domain $\mathcal{D}^{\hspace{0.02cm}\scriptscriptstyle \beta}_{\scriptscriptstyle N}\hspace{0.03cm}$, and then writing down the corresponding version of Lemma \ref{lem:rgdsym}, which directly yields Proposition \ref{prop:acosh}. Fortunately, all the necessary information can be found in~\cite{piatetski}\cite{helgason}. \\[0.1cm]
\textbf{The symmetric pair\,:} $\mathcal{D}^{\hspace{0.02cm}\scriptscriptstyle \beta}_{\scriptscriptstyle N}$ is associated to the symmetric pair $(G,K)$, where $G$ and $K$ are groups of $2N \times 2N$ complex matrices $g$, defined in the following way (here, $^*$ denotes the complex conjugate). 
$$
\begin{array}{rl}
G_{\beta = 2} = & \left\lbrace g = \left(\begin{array}{cc} A & B \\ C & D\end{array}\right): g\hspace{0.02cm}\mathrm{P}g^\dagger = \mathrm{P} \right\rbrace \\[1cm]
G_{\beta = 1} = & \left\lbrace g = \left(\begin{array}{cc} A & B \\ C & D\end{array}\right): \begin{array}{r} g\hspace{0.02cm}\mathrm{P}g^\dagger = \mathrm{P} \\ \text{ and } g\hspace{0.02cm} \mathrm{S}g^T = \mathrm{S}\end{array}\right\rbrace
\end{array}
$$

$$
\begin{array}{rl}
K_{\beta = 2} = & \left\lbrace k = \left(\begin{array}{cc} U & \\ & V \end{array}\right):U\hspace{0.02cm},V \in U(N) \right\rbrace
 \\[1cm]
K_{\beta = 1} = & \left\lbrace k = \left(\begin{array}{cc} U & \\ & U^* \end{array}\right):U\in U(N) \right\rbrace 
\end{array}
$$
\vspace{0.2cm}

\noindent where $\mathrm{P}$ and $\mathrm{S}$ denote the following matrices

$$
\mathrm{P} = \left(\begin{array}{cc} \mathrm{I}_{\scriptscriptstyle N} & \mathrm{0}_{\scriptscriptstyle N} \\ \mathrm{0}_{\scriptscriptstyle N} & -\mathrm{I}_{\scriptscriptstyle N} \end{array}\right) \,;\,
\mathrm{S} = 
\left(\begin{array}{cc} \mathrm{0}_{\scriptscriptstyle N} & \mathrm{I}_{\scriptscriptstyle N} \\ -\mathrm{I}_{\scriptscriptstyle N} & \mathrm{0}_{\scriptscriptstyle N} \end{array}\right)
$$
\vspace{0.1cm}

\noindent These groups $G$ and $K$ act on $\mathcal{D}^{\hspace{0.02cm}\scriptscriptstyle \beta}_{\scriptscriptstyle N}$ by matrix fractional transformations
$$
g \cdot \Omega = (A\Omega + B)(C\Omega + D)^{-1} \hspace{0.2cm} g\in G \text{ and  }\Omega \in \mathcal{D}^{\hspace{0.02cm}\scriptscriptstyle \beta}_{\scriptscriptstyle N}
$$
\noindent \textbf{The Lie algebras\,:} these are given by
$$
\begin{array}{rl}\mathfrak{g}_{\hspace{0.02cm}\beta = 2} = &\lbrace (\gamma,\delta,\varepsilon): \gamma + \gamma^\dagger = \varepsilon + \varepsilon^\dagger = \mathrm{0}_{\scriptscriptstyle N}\rbrace \\[0.2cm]
\mathfrak{g}_{\hspace{0.02cm}\beta = 1} = & \left\lbrace (\gamma,\delta,\varepsilon): \begin{array}{rl}\gamma + \gamma^\dagger = \varepsilon + \varepsilon^\dagger & = \mathrm{0}_{\scriptscriptstyle N} \\ \delta - \delta^T = \gamma + \varepsilon^T  & = \mathrm{0}_{\scriptscriptstyle N}\end{array} \right\rbrace \!\!
\end{array}
$$
$$
\begin{array}{rl}
\mathfrak{k}_{\hspace{0.02cm}\beta = 2} = & \lbrace (\gamma,\mathrm{0}_{\scriptscriptstyle N},\varepsilon): \gamma + \gamma^\dagger = \varepsilon + \varepsilon^\dagger = \mathrm{0}_{\scriptscriptstyle N}\rbrace \\ \mathfrak{k}_{\hspace{0.02cm}\beta = 1} = & \lbrace (\gamma,\mathrm{0}_{\scriptscriptstyle N},\gamma^*): \gamma + \gamma^\dagger = \mathrm{0}_{\scriptscriptstyle N}\rbrace \end{array}
$$
in terms of the notation
$$
(\gamma,\delta,\varepsilon) = \left(\begin{array}{cc} \gamma & \delta \\ \delta^\dagger & \varepsilon \end{array}\right)
$$

\noindent \textbf{The subspaces $\mathfrak{p}$ and $\mathfrak{a}$\,:} these are given by
$$
\begin{array}{rl}
\mathfrak{p}_{\hspace{0.02cm}\beta = 2} = & \lbrace (\mathrm{0}_{\scriptscriptstyle N},\delta,\mathrm{0}_{\scriptscriptstyle N}) \rbrace \\
\mathfrak{p}_{\hspace{0.02cm}\beta = 1} = & \lbrace (\mathrm{0}_{\scriptscriptstyle N},\delta,\mathrm{0}_{\scriptscriptstyle N}):  \delta - \delta^T = \mathrm{0}_{\scriptscriptstyle N}\rbrace
\end{array}
$$
$$
\mathfrak{a}_{\hspace{0.02cm}\beta = 2} = 
\mathfrak{a}_{\hspace{0.02cm}\beta = 1} = \lbrace a =  (\mathrm{0}_{\scriptscriptstyle N},r,\mathrm{0}_{\scriptscriptstyle N}): r \text{ real diagonal}\rbrace \hspace{1.5cm}\,
$$

\noindent \textbf{Lemma \ref{lem:rgdsym}\,:} the only information still needed is the positive roots. Matrix multiplication shows that for any $a \in \mathfrak{a}$, $\mathrm{ad}_a$ has the following eigenvalues and eigenvectors in $\mathfrak{g}_{\hspace{0.02cm}\beta = 2}\hspace{0.03cm}$.
$$
\begin{array}{lcl}
\xi = \left(\omega_{ij},\tau_{ij},\omega_{ij}\right) & \Longrightarrow & \mathrm{ad}_a(\xi) = (r_i - r_j)\hspace{0.02cm}\xi\\
\xi = \mathrm{i}\left(\tau_{ij},\omega_{ij},\tau_{ij}\right) & \Longrightarrow & \mathrm{ad}_a(\xi) = (r_i - r_j)\hspace{0.02cm}\xi
\\[0.1cm]
\xi = \mathrm{i}\left(-\tau_{ij},\tau_{ij},\tau_{ij}\right) & \Longrightarrow & \mathrm{ad}_a(\xi) = (r_i + r_j)\hspace{0.02cm}\xi \\ 
\xi = \left(\omega_{ji},\omega_{ij},\omega_{ij}\right) & \Longrightarrow & \mathrm{ad}_a(\xi) = (r_i + r_j)\hspace{0.02cm}\xi
\end{array}
$$
where $\tau_{ij} = e_{ij} + e_{ji}$ and $\omega_{ij} = e_{ij} - e_{ji}\hspace{0.03cm}$, with $e_{ij}$ a matrix all of whose entries are zero, except the one on line $i$ and column $j$, which is equal to $1$, and where $r = \mathrm{diag}(r_{\scriptscriptstyle 1},\ldots,r_{\scriptscriptstyle N})$. 

This shows that the positive roots are $\rho(a) = r_i - r_j$ where $i < j$ and $\rho(a) = r_i + r_j$ where $i \leq j$, which all have multiplicity $m_\rho = 2 $ (this is $\beta$), except for $\rho(a) = 2r_i$ (this is $r_i + r_j$ when $i = j$), which has multiplicity $m_\rho = 1$. Moreover, of the above eigenvectors, only the ones in the left column belong to $\mathfrak{g}_{\hspace{0.02cm}\beta = 1}\hspace{0.03cm}$. Thus, in the $\beta = 1$ case, all the multiplicities $m_\rho$ are equal to $1$. Finally, note that, from the power series of the matrix exponential,
\begin{equation} 
\exp(a) = \left(\begin{array}{cc} \cosh(r) & \sinh(r) \\ \sinh(r) & \cosh(r) \end{array}\right) \,;\, a = (\mathrm{0}_{\scriptscriptstyle N},r,\mathrm{0}_{\scriptscriptstyle N})
\end{equation}

\noindent From the above form of the action of $G$ on $\mathcal{D}^{\hspace{0.02cm}\scriptscriptstyle \beta}_{\scriptscriptstyle N\hspace{0.02cm}}$, it is then straightforward that the representation (\ref{eq:polarcoordinates}) is the same as (\ref{eq:takagi}). This implies the first part of Proposition \ref{prop:acosh} ($U, V$ are uniformly distributed on $U(N)$). For the second part, it follows from (\ref{eq:aadensity}) that, in the present case,
\begin{align}
\nonumber p(r|\sigma) \propto \prod_{i} \exp\left[-\frac{r^{\scriptscriptstyle 2}_i}{2\sigma^{\scriptscriptstyle 2}}\right]\sinh(2r_i)\\
 \prod_{i<j} \left(\sinh\left|r_i-r_j\right|\sinh\left|r_i+r_j\right|\right)^\beta
\end{align}
This yields (\ref{eq:acoshnormal2}), by using $2\sinh\left|r_i-r_j\right|\sinh\left|r_i+r_j\right| = |\cosh(2r_i) - \cosh(2r_j)|$, and introducing the change of variables $x_i = \cosh(2r_i)$. 
    
\vfill
\pagebreak

\section*{Acknowledgment}
S.H. is supported by the Science and Technology Facilities Council (STFC) and St. John's College, Cambridge. S.H. also benefited from partial support from the Cambridge Mathematics Placement (CMP) Programme and the European Research Council under the Advanced ERC Grant Agreement Switchlet n.670645. C.M.
was supported by a Henslow Fellowship from the Cambridge Philosophical Society as well as a Presidential Postdoctoral Fellowship at NTU and is grateful for support from Fitzwilliam College, Cambridge.

\bibliographystyle{IEEEtran}
\bibliography{References}

\begin{thebibliography}{10}
\providecommand{\url}[1]{#1}
\csname url@samestyle\endcsname
\providecommand{\newblock}{\relax}
\providecommand{\bibinfo}[2]{#2}
\providecommand{\BIBentrySTDinterwordspacing}{\spaceskip=0pt\relax}
\providecommand{\BIBentryALTinterwordstretchfactor}{4}
\providecommand{\BIBentryALTinterwordspacing}{\spaceskip=\fontdimen2\font plus
\BIBentryALTinterwordstretchfactor\fontdimen3\font minus
  \fontdimen4\font\relax}
\providecommand{\BIBforeignlanguage}[2]{{%
\expandafter\ifx\csname l@#1\endcsname\relax
\typeout{** WARNING: IEEEtran.bst: No hyphenation pattern has been}%
\typeout{** loaded for the language `#1'. Using the pattern for}%
\typeout{** the default language instead.}%
\else
\language=\csname l@#1\endcsname
\fi
#2}}
\providecommand{\BIBdecl}{\relax}
\BIBdecl

\bibitem{Cheng2013}
\BIBentryALTinterwordspacing
G.~Cheng and B.~C. Vemuri, ``A novel dynamic system in the space of {SPD}
  matrices with applications to appearance tracking,'' \emph{SIAM Journal on
  Imaging Sciences}, vol.~6, no.~1, pp. 592--615, 2013. [Online]. Available:
  \url{https://doi.org/10.1137/110853376}
\BIBentrySTDinterwordspacing

\bibitem{Sa16}
S.~Said, L.~Bombrun, Y.~Berthoumieu, and J.~H. Manton, ``Riemannian {Gaussian}
  distributions on the space of symmetric positive definite matrices,''
  \emph{IEEE Transactions on Information Theory}, vol.~63, no.~4, pp.
  2153--2170, 2017.

\bibitem{Sa17}
S.~Said, H.~Hajri, L.~Bombrun, and B.~C. Vemuri, ``Gaussian distributions on
  {Riemannian} symmetric spaces: Statistical learning with structured
  covariance matrices,'' \emph{IEEE Transactions on Information Theory},
  vol.~64, no.~2, pp. 752--772, 2018.

\bibitem{zanini2018}
P.~Zanini, M.~Congedo, C.~Jutten, S.~Said, and Y.~Berthoumieu, ``Transfer
  learning: {A} {Riemannian} geometry framework with applications to
  brain–computer interfaces,'' \emph{IEEE Transactions on Biomedical
  Engineering}, vol.~65, no.~5, pp. 1107--1116, 2018.

\bibitem{mathieu2019}
E.~Mathieu, C.~Le~Lan, C.~J. Maddison, R.~Tomioka, and Y.~W. Teh, ``Continuous
  hierarchical representations with {P}oincar\'{e} variational auto-encoders,''
  in \emph{Advances in Neural Information Processing Systems}, vol.~32.\hskip
  1em plus 0.5em minus 0.4em\relax Curran Associates, Inc., 2019.

\bibitem{SH21}
S.~Heuveline, S.~Said, and C.~Mostajeran, ``Gaussian distributions on
  {R}iemannian symmetric spaces in the large ${N}$ limit,'' in \emph{Geometric
  Science of Information}, F.~Nielsen and F.~Barbaresco, Eds.\hskip 1em plus
  0.5em minus 0.4em\relax Cham: Springer International Publishing, 2021, pp.
  20--28.

\bibitem{Ma05}
M.~Mari\~no, \emph{Chern-Simons theory, matrix models, and topological
  strings}.\hskip 1em plus 0.5em minus 0.4em\relax Oxford University Press,
  2005.

\bibitem{Ti20}
\BIBentryALTinterwordspacing
L.~Santilli and M.~Tierz, ``Riemannian {Gaussian} distributions, random matrix
  ensembles and diffusion kernels,'' \emph{Nuclear Physics B}, vol. 973, p.
  115582, 2021. [Online]. Available:
  \url{https://www.sciencedirect.com/science/article/pii/S0550321321002790}
\BIBentrySTDinterwordspacing

\bibitem{Forrester1989}
\BIBentryALTinterwordspacing
P.~J. Forrester, ``Vicious random walkers in the limit of a large number of
  walkers,'' \emph{Journal of Statistical Physics}, vol.~56, no.~5, pp.
  767--782, 1989. [Online]. Available: \url{https://doi.org/10.1007/BF01016779}
\BIBentrySTDinterwordspacing

\bibitem{Forrester2021}
\BIBentryALTinterwordspacing
------, ``Global and local scaling limits for the $\beta$ = 2
  {Stieltjes}–{Wigert} random matrix ensemble,'' \emph{Random Matrices:
  Theory and Applications}, vol.~11, no.~02, 2022. [Online]. Available:
  \url{https://doi.org/10.1142/S2010326322500204}
\BIBentrySTDinterwordspacing

\bibitem{amaribook}
S.~I. Amari, \emph{Information geometry and its applications}.\hskip 1em plus
  0.5em minus 0.4em\relax Springer, 2016.

\bibitem{Pennec2006}
\BIBentryALTinterwordspacing
X.~Pennec, ``Intrinsic statistics on {Riemannian} manifolds: Basic tools for
  geometric measurements,'' \emph{Journal of Mathematical Imaging and Vision},
  vol.~25, no.~1, pp. 127--154, 2006. [Online]. Available:
  \url{https://doi.org/10.1007/s10851-006-6228-4}
\BIBentrySTDinterwordspacing

\bibitem{Pennec2019}
X.~Pennec, S.~Sommer, and T.~Fletcher, \emph{Riemannian geometric statistics in
  medical image analysis}.\hskip 1em plus 0.5em minus 0.4em\relax Academic
  Press, 2020.

\bibitem{edelman}
\BIBentryALTinterwordspacing
A.~Edelman and S.~Jeong, ``Fifty three matrix factorizations: A systematic
  approach,'' 2021. [Online]. Available: \url{https://arxiv.org/abs/2104.08669}
\BIBentrySTDinterwordspacing

\bibitem{meckes}
E.~S. Meckes, \emph{The Random Matrix Theory of the Classical Compact Groups},
  ser. Cambridge Tracts in Mathematics.\hskip 1em plus 0.5em minus 0.4em\relax
  Cambridge University Press, 2019.

\bibitem{Sommer2015}
S.~Sommer, ``Anisotropic distributions on manifolds: Template estimation and
  most probable paths,'' in \emph{Information Processing in Medical Imaging},
  S.~Ourselin, D.~C. Alexander, C.-F. Westin, and M.~J. Cardoso, Eds.\hskip 1em
  plus 0.5em minus 0.4em\relax Cham: Springer International Publishing, 2015,
  pp. 193--204.

\bibitem{deift}
P.~Deift, \emph{Orthogonal polynomials and random matrices: a
  {Riemann}-{Hilbert} approach}.\hskip 1em plus 0.5em minus 0.4em\relax
  American Mathematical Society, 1999.

\bibitem{mtnshmm}
\BIBentryALTinterwordspacing
S.~Said, N.~L. Bihan, and J.~H. Manton, ``Hidden {Markov} chains and fields
  with observations in {Riemannian} manifolds,'' vol.~54, no.~9, 2021, pp.
  719--724, 24th International Symposium on Mathematical Theory of Networks and
  Systems MTNS 2020. [Online]. Available:
  \url{https://www.sciencedirect.com/science/article/pii/S2405896321006200}
\BIBentrySTDinterwordspacing

\bibitem{QT21}
Q.~Tupker, S.~Said, and C.~Mostajeran, ``Online learning of {R}iemannian hidden
  {M}arkov models in homogeneous {H}adamard spaces,'' in \emph{Geometric
  Science of Information}, F.~Nielsen and F.~Barbaresco, Eds.\hskip 1em plus
  0.5em minus 0.4em\relax Cham: Springer International Publishing, 2021, pp.
  37--44.

\bibitem{habilitation}
S.~Said, ``Statistical models and probabilistic methods on {Riemannian}
  manifolds,'' 2021, {Universit\'e de Bordeaux}, {Habilitation} thesis (HDR).

\bibitem{jeuris1}
\BIBentryALTinterwordspacing
B.~Jeuris and R.~Vandebril, ``The {K\"ahler} mean of block {Toeplitz} matrices
  with {Toeplitz} structured blocks,'' \emph{SIAM Journal on Matrix Analysis
  and Applications}, vol.~37, no.~3, pp. 1151--1175, 2016. [Online]. Available:
  \url{https://doi.org/10.1137/15M102112X}
\BIBentrySTDinterwordspacing

\bibitem{yannthesis}
Y.~Cabanes, ``Multidimensional complex stationary centred {Gaussian}
  autoregressive model classification\,: Applications for audio and radar
  machine learning in hyperbolic and {Siegel} spaces,'' Ph.D. dissertation,
  Institut de Math\'ematiques de Bordeaux, 2022, {M.} Arnaudon, F. Barbaresco,
  J. Bigot (supervisors).

\bibitem{piatetski}
I.~I. Piatetskii-Shapiro, \emph{Automorphic functions and the geometry of
  classical domains}.\hskip 1em plus 0.5em minus 0.4em\relax Gordon and Breach
  Science Publishers, 1969.

\bibitem{helgason}
S.~Helgason, \emph{Differential Geometry and Symmetric Spaces}.\hskip 1em plus
  0.5em minus 0.4em\relax Academic Press, 1962.

\end{thebibliography}

\begin{IEEEbiographynophoto}{Salem Said}
has been a CNRS research scientist since 2014. He defended his Habilitation thesis, at the university of Bordeaux, in 2021, and moved to the university of Grenoble, in 2022, where he is now based at Laboratoire Jean Kuntzmann. His research focuses on statistical inference and stochastic optimisation on Riemannian manifolds, and on applications to brain-computer interface analysis.
\end{IEEEbiographynophoto}
\begin{IEEEbiographynophoto}{Simon Heuveline}
is a PhD student under the supervision of David Skinner in the Department of Applied Mathematics and Theoretical Physics (DAMTP) at the University of Cambridge. He is a member of St. John's College and his research focuses on topics surrounding string theory, flat space holography and twistor theory. He studied mathematics and physics as an undergraduate at the University of Heidelberg and the University of Cambridge before doing graduate work in differential geometry and mathematical physics. 
\end{IEEEbiographynophoto}
\begin{IEEEbiographynophoto}{Cyrus Mostajeran}
is a Presidential Postdoctoral Fellow at the School of Physical and Mathematical Sciences at Nanyang Technological University (NTU) in Singapore. He studied mathematics as an undergraduate at Balliol College in the University of Oxford before doing graduate work in mathematics, physics, and engineering, earning a PhD in Information Engineering from the University of Cambridge in 2018. He was an Early Career Research Fellow at Fitzwilliam College in the University of Cambridge between 2018 and 2022.  A common theme of his research is the application of differential geometry to problems arising in statistics, optimization, materials science, and robotics.
\end{IEEEbiographynophoto}

\end{document}